\documentclass[11pt,leqno]{article}
\usepackage[centertags]{amsmath}
\usepackage{amsfonts}
\usepackage{amsmath}
\usepackage{amssymb}
\usepackage{latexsym}
\usepackage{amsthm}
\usepackage{mathrsfs}
\usepackage{newlfont}
\usepackage[latin1]{inputenc}
\usepackage{bm}

\usepackage{epsfig}
\usepackage{graphicx}

\newtheorem{teor}{Theorem}[section]

\newtheorem{defi}{Definition}[section]
\newtheorem{lema}{Lemma}[section]

\newtheorem{nota}{Remark}[section]

\newcommand{\LWd}{L_{\bm W}^2(S,\mathbb{C}^{N\times N})}

\numberwithin{equation}{section} %%%%%%numera las ecuaciones (1.1),
\title{Spectral methods for bivariate Markov processes with diffusion and discrete components and a variant of the Wright-Fisher model\thanks{The work of the author is partially supported by the Ministry of Science and Innovation of Spain, grant MTM2009-12740-C03-02, by Junta de Andaluc\'{i}a grants  FQM-262, P06-FQM-01735 and P09-FQM-4643, and by Subprograma de estancias de movilidad posdoctoral en el extranjero, MICINN, ref. -2008-0207.\newline \textbf{AMS Subject
Classifications.} 60K37, 60J60, 42C05. \newline \textbf{Key words.} Bivariate Markov processes, switching diffusions, matrix-valued orthogonal functions, Wright-Fisher models.}}

\author{Manuel D. de la Iglesia \\ \footnotesize Departamento de An\'{a}lisis Matem\'{a}tico. Universidad de Sevilla \\ \footnotesize\ Apdo (P. O. BOX) 1160. 41080 Sevilla. Spain ({\tt mdi29@us.es}).}

\begin{document}

\maketitle

\begin{abstract}
The aim of this paper is to study differential and spectral properties of the infinitesimal operator of two dimensional Markov processes with diffusion and discrete components. The infinitesimal operator is now a second-order differential operator with matrix-valued coefficients, from which we can derive backward and forward equations, a spectral representation of the probability density, study recurrence of the process and the corresponding invariant distribution. All these results are applied to an example coming from group representation theory which can be viewed as a variant of the Wright-Fisher model involving only mutation effects.
\end{abstract}

\section{Introduction}\label{INTRO}

The connection between stochastic processes and orthogonal polynomials is very well known. For a very good account see \cite{Sch}. Among them we can find diffusion processes where the state space is an interval of the real line, i.e. $\mathcal{S}\subset\mathbb{R}$, and the parameter set is $\mathcal{T}=[0,\infty)$, i.e. continuous time. In this case, the corresponding infinitesimal operator $\mathcal{A}$ comes in terms of a second-order differential operator. It is possible to calculate an expression of the transition probability density by \emph{spectral methods}, i.e. in terms of the eigenfunctions of the infinitesimal operator $\mathcal{A}$ (see Section 8, Chapter V of \cite{BW} or Section 13, Chapter 15 of \cite{KT}), among other properties. In some cases, the invariant distribution can also be given in terms of the associated orthogonality measure of the eigenfunctions. Prominent examples connected with spectral methods are the Orstein-Uhlenbeck process, population growth models or Wright-Fisher models.

Continuous-time bivariate Markov processes with diffusion and discrete components are a natural extension, where now the state space is two dimensional, i.e. $\mathcal{S}=S\times\{1,2,\ldots,N\}$, where $S\subset\mathbb{R}$ and $N$ is a positive integer. The first component is continuous, which we will call the \emph{diffusion} component, while the second one is discrete and which we will call the \emph{phase}. We will denote these bivariate Markov processes by $(X_t, Y_t), t\geq0$. At first it would seem that this hybrid process is just blending a regular scalar diffusion process and a continuous time Markov process, but there are more subtleties specially in the case when both components are dependent each other.

The ideas behind these coupled stochastic processes may be traced back to what is known in the literature as \emph{random evolutions}. These were introduced by R. Griego and R. Hersh in 1969 \cite{GH1}, although the seminal work was given by M. Kac several years before (see \cite{He2} for an historical account). The main idea of a random evolution is a model for a dynamical system whose equation of state is subject to random variation. Other authors like G. Papanicolaou, M. Pinsky, T. Kurtz or R. Varadhan, to mention a few, were involved in this field during the 1960s and the 1970s. For a very good account see the Hersh survey \cite{He1} and the books \cite{Pap, Pins}. In the 1980s, the Kiev school led by Korolyuk and Swishchuk put these processes in the framework of semi-Markov processes (see \cite{Swi} for a recent book on this subject).

When the equation of state is governed by diffusion processes, we are in the case of bivariate Markov processes, the objects we will be interested in this paper. Now the infinitesimal operator $\mathcal{A}$ comes in terms of a second-order differential operator with matrix-valued coefficients. The basic definitions and properties of these bivariate Markov processes were introduced in \cite{B} (see also Section \ref{PRE} below). Here the author focused on the calculation of the first-passage time distribution out of the phases using the resolvent of the corresponding infinitesimal operator. This model was shown to be applicable to describing the progression of HIV in an infected individual. These processes have received other names in the literature, like \emph{diffusions in random environments} or \emph{switching diffusion models} (see for instance \cite{GHO}).

%Also in \cite{Bl} the author describes an approach to Bayesian inference for these models and applies it to study data arising in certain radio-tracking experiments. These processes are also related with diffusions in random environments and switching diffusion models (see \cite{Bl} and references therein). The difference with the work in \cite{B} is that the discrete component of the process acts independently of the continuous component.

It is not difficult to foresee the applications that these bivariate Markov processes may have, as in random evolutions or switching diffusion models, in fields like population dynamics, insurance or financial mathematics.

%If the first component of these bivariate Markov processes is taken to be the set of nonnegative integers, then they are called quasi-birth-and-death processes. They have been connected with queuing problems in network theory as well as the general field of communication systems (see \cite{LR1, N}). The approach of these processes using spectral methods was initially raised independently by \cite{DRSZ} and \cite{G2} in discrete time (for continuous time see \cite{DR}).

In this paper we will focus on differential and spectral properties of the infinitesimal operator generated by a bivariate Markov process with diffusion and discrete components. Surprisingly enough, this has not been explored in depth so far. The aim of this paper is to fill this gap using spectral analysis of differential operators with matrix-valued coefficients.% (see \cite{DG1, GPT1, GPT5}).

In the first part, after some preliminary considerations in Section \ref{PRE} (already introduced in \cite{B}), we will derive, in Section \ref{NEW}, the corresponding backward and forward equations, a representation of the (matrix-valued) probability density in terms of the eigenfunctions of the infinitesimal operator $\mathcal{A}$ and some differential equations associated with some functionals, including the invariant distribution of these processes.

In the second part, in Section \ref{EXAM}, we will study an example (coming from group representation theory \cite{GPT1, PT1}) of these bivariate Markov processes which can be viewed as a variant of the Wright-Fisher model involving only mutation effects. We will use the results in Sections \ref{PRE} and \ref{NEW} to study the probabilistic aspects of this example, specially the spectral representation of the probability density and an explicit expression of the invariant distribution.

\section{Preliminaries}\label{PRE}

Let $(X_t, Y_t), t\geq0$ be a time-homogeneous Markov process assuming values in $S\times\{1,2,\ldots,N\}$, where $S\subset\mathbb{R}$ and $N$ is a positive integer. We use preferably boldface capital letters to denote matrices, boldface lowercase letters to denote vectors and standard font for scalars. We define the transition probability density as a matrix-valued function $\bm P(t;x,A)=(P_{ij}(t;x,A))$, defined for every $t\geq 0, x\in S$, and real Borel set $A$, as a matrix whose entry $(i,j)$ is given by
\begin{equation}\label{TPD}
 P_{ij}(t;x,A)=\mbox{Pr}\{X_t\in A, Y_t=j|X_0=x, Y_0=i\}.
\end{equation}
For simplicity, we will assume that every entry $P_{ij}(t;x,dy)$ has a density, also called $P_{ij}(t;x,y)$, in such a way that we can write
$$
 P_{ij}(t;x,A)=\int_AP_{ij}(t;x,y)dy.
$$

Observe that $\bm P(t;x,y)$ is not a probability density in the common sense, but we have that all entries are positive and for every real Borel set $A\subset S$
$$
\left(\int_A\bm P(t;x,y)dy\right)\bm e_N\leq\bm e_N,
$$
where $\bm e_N$ is $N$-dimensional column vector with all entries equal to 1. The equality holds when $A=S$.
% In particular we have for every row $i$ that
%$$
%\sum_{j=1}^N P_{ij}(t;x,S)=1, \quad t\geq0,
%$$
%that is, starting at any $X_0=x, Y_0=i$, the probability of finding the process $X_t$ at any other real number is the sum of the probabilities that $Y_t$ is at one of the discrete phases, which is 1.

It is assumed that $\bm{P}(t;x,y)$ satisfies the Chapman-Kolmogorov equation
\begin{equation}\label{CKEQ}
 \bm{P}(s+t;x,y)=\int_S\bm{P}(s;x,z)\bm{P}(t;z,y)dz.
\end{equation}

In a similar way as the classical theory of diffusion processes we have to make the following assumptions:

\begin{enumerate}
 \item For every $\epsilon>0$
\begin{equation}\label{AS1}
 \lim_{h\rightarrow 0}\frac{1}{h}\bm{P}(h;x,\{y: |y-x|>\epsilon\})=\bm 0.
\end{equation}
\item There exist $N\times N$ matrix-valued functions $\bm{A}(x)$ (with positive entries) and $\bm{B}(x)$, $x\in S$, such that for any $\epsilon>0$
\begin{equation}\label{AS21}
 \lim_{h\rightarrow 0}\frac{1}{h}\left[\int_{\{y: |y-x|\leq\epsilon\}}(y-x)\bm{P}(h;x,y)dy\right]=\bm{B}(x)
\end{equation}
and
\begin{equation}\label{AS22}
 \lim_{h\rightarrow 0}\frac{1}{h}\left[\int_{\{y: |y-x|\leq\epsilon\}}(y-x)^2\bm{P}(h;x,y)dy\right]=\bm{A}(x).
\end{equation}
\item There exists an $N\times N$ matrix-valued function $\bm{Q}(x)=(Q_{ij}(x))$ with $Q_{ij}$ real-valued such that
\begin{equation}\label{AS3}
 \lim_{h\rightarrow 0}\frac{1}{h}\left[\bm{P}(h;x,S)-\bm I\right]=\bm{Q}(x).
\end{equation}
%\item For every $\epsilon>0$ and integer $p\geq3$
%\begin{equation}\label{AS4}
% \lim_{h\rightarrow 0}\frac{1}{h}\left[\int_{\{y: |y-x|\leq\epsilon\}}|y-x|^p\bm{P}(h;x,dy)\right]=0.
%\end{equation}
\end{enumerate}

Observe that (\ref{AS1}) gives continuous sample paths for the process $X_t$ while (\ref{AS21}) and (\ref{AS22}) are matrix-valued versions of the corresponding hypothesis employed in diffusion processes (see, for example, \cite{BW, KT}). Note that (\ref{AS3}) is a more general form of a corresponding hypothesis employed in continuous time Markov chains. The formulation here is different because $\bm Q$ \emph{depends on} $x$. In particular we have that, for fixed $x$,
\begin{eqnarray*}
Q_{ij}(x) &\leq&0\quad\mbox{for}\quad i=j, \\
Q_{ij}(x) &\geq& 0\quad\mbox{for}\quad i\neq j;
\end{eqnarray*}
and for each $i$
\begin{equation*}\label{sum0}
 -Q_{ii}(x)=\sum_{i\neq j}Q_{ij}(x),
\end{equation*}
i.e., every off diagonal entry of $\bm{Q}(x)$ is a positive function and the sum of every row of $\bm{Q}(x)$ is 0, for every $x$, i.e., $\bm{Q}(x)\bm e_N=\bm0$.

\medskip

We will denote by $\mathfrak{B}(S^N)$ the set of all column vector-valued functions where all components are real-valued, bounded and Borel measurable functions on $S$.

For any $\bm{f}(x)\in\mathfrak{B}(S^{N})$ we define the transition operator as
\begin{equation}\label{MVE}
\bm{T}_t \bm{f}(x)=\bm{E}[\bm{f}(X_t)|X_0=x]\doteq\int_S\bm{P}(t;x,y)\bm{f}(y)dy.
\end{equation}
%Observe that every entry $i$ of (\ref{MVE}) gives
%$$
%E[f_{Y_t}(X_t)|X_0=x,Y_0=i]=\sum_{k=1}^N\int_SP_{ik}(t;x,y)f_{k}(y)dy.
%$$

The transition operator (\ref{MVE}) can be applied to matrix-valued functions $\bm{F}(x)\in\mathfrak{B}(S^{N\times N})$, with similar properties as $\mathfrak{B}(S^N)$. We just have to perform the operator column by column. For simplicity we will use vector-valued analysis, but sometimes it will be more convenient to use a matrix-valued approach. Both are vector spaces with a norm over $\mathbb{R}$.

The family of transition operators $\{\bm{T}_t: t>0\}$ has the semigroup property
$$
\bm{T}_{s+t}=\bm{T}_s\bm{T}_t,
$$
as a consequence of the Chapman-Kolmogorov equation \eqref{CKEQ}. In particular, this implies that the transition operators commute. Therefore, the behavior of $\bm{T}_t$ near $t=0$ completely determines the semigroup, as in the scalar situation.

The infinitesimal generator $\mathcal{A}$ of the process is defined for each smooth vector-valued function $\bm{f}(x)$ as
\begin{equation}\label{infgen}
 (\mathcal{A}\bm{f})(x)=\lim_{h\rightarrow 0}\frac{1}{h}\left[\bm{T}_t \bm{f}(x)-\bm{f}(x)\right]=\lim_{h\rightarrow 0}\frac{1}{h}\left[\int_S\bm{P}(h;x,y)\bm{f}(y)dy-\bm{f}(x)\right].
\end{equation}

The following two important results were proved in \cite{B}:
\begin{lema}[Lemma 2.1 of \cite{B}]
 $\bm{A}(x)$ and $\bm{B}(x)$, $x\in S$, are diagonal matrices.
\end{lema}
\begin{teor}[Theorem 2.1 of \cite{B}]
 For any vector-valued function $\bm{f}(x), x\in S,$ whose entries are real-valued, bounded and have continuous second derivatives, the infinitesimal operator $\mathcal{A}$ defined in (\ref{infgen}) is of the form
\begin{equation}\label{ForIO}
  (\mathcal{A}\bm{f})(x)=\frac{1}{2}\bm{A}(x)\bm{f}''(x)+\bm{B}(x)\bm{f}'(x)+\bm{Q}(x)\bm{f}(x).
\end{equation}
\end{teor}

For simplicity, since the coefficients $\bm{A}(x)$ and $\bm{B}(x)$ are diagonal matrices, we will denote the corresponding diagonal entries by $\sigma_i^2(x)$ and $\tau_i(x)$, $i=1,\ldots,N$, respectively.

\medskip

\begin{figure}[h]
\begin{center}
\vspace{-6.0cm}
\includegraphics[height=20cm]{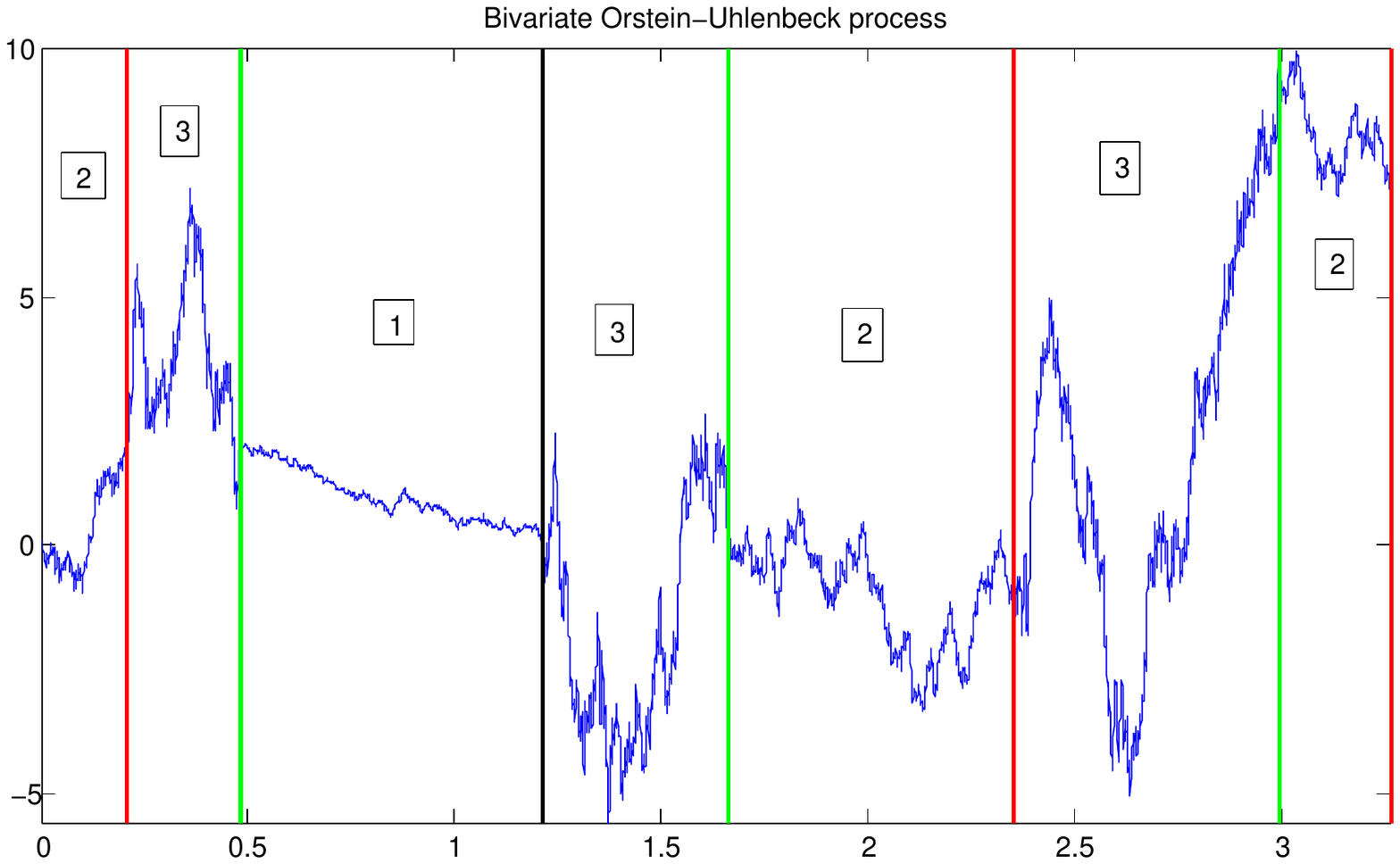}
\vspace{-7.0cm}
\end{center}
\caption{Sample path of a bivariate Markov process with 3 phases.}
\label{OU3}
\end{figure}

Now we will give a stochastic representation of these processes, already given in Section 3 of \cite{B}. We refer to Figure \ref{OU3} for a better comprehension of how these processes work. 

Suppose the process $(X_t,Y_t)$ starts at state $X_0=x$ and phase $Y_0=i$. Let $X^i_t$ be the diffusion process starting at the point $x$ and evolving according to the infinitesimal generator
\begin{equation}\label{dinfgen}
 \frac{1}{2}\sigma_i^2(x)\frac{d^2}{dx^2}+\tau_i(x)\frac{d}{dx}+Q_{ii}(x)\frac{d^0}{dx^0}.
\end{equation}
$X^i_t$ is a diffusion process with diffusion and drift coefficients $\sigma_i^2(x)$ and $\tau_i(x)$, respectively, and a \emph{killing} coefficient $-Q_{ii}(x)\geq0$. Since the diffusion $X^i_t$ has a lifetime $t_i$ (chosen according to exponential holding time of $-Q_{ii}(X_{t})$, see pp. 314 of \cite{KT}), the process is defined on $[0,t_i]$ as $X_t=X^i_t, Y_t=i$ (if $t_i=\infty$, then $X_t=X^i_t$ and $Y_t=i$ for all $t>0$). At time $t=t_i$, $Y_t$ changes from phase $i$ to another phase $j\neq i$. If $-Q_{ii}(X_{t_i})>0$, there is a probability distribution over the phases $j\neq i$, defined as $-Q_{ij}(X_{t_i})/Q_{ii}(X_{t_i}), j\neq i$, and $Y_{t}$ is selected in accordance with this distribution.

At the moment $t_i$, the diffusion component $X_t$ is changed from $X^i_t$ to $X^j_t$, starting at $X^i_{t_i}$ and evolving according with the generator (\ref{dinfgen}) with index $j$ in the place of $i$. The process is killed at the time $t_i+t_j$, where $t_j$ is the lifetime of $X^j_t$ and the process is defined on $[t_i,t_i+t_j]$ as $X_t=X^j_t, Y_t=j$. If $-Q_{jj}(X_{t_i+t_j})>0$, there is a probability distribution over the phases $k\neq j$, defined as $-Q_{jk}(X_{t_i+t_j})/Q_{jj}(X_{t_i+t_j}), k\neq j$, and $Y_t$ is selected in accordance with this distribution. Then the diffusion changes from $X^j_t$ to $X^k_t$ at time $t_i+t_j$, where $X^k_t$ starts at $X_{t_i+t_j}$ and has generator (\ref{dinfgen}) with $k$ in the place of $i$, and so on.

Observe that the sample path of the process is defined as a piecewise function and is always \emph{continuous} for $t>0$. The only difference with respect to a regular diffusion process is that this process can move through $N$ different phases with different diffusion and drift coefficients in general. The transitions through all the phases and the lifetime spent at each phase, which depend on the position $X_t$, are derived in terms of the coefficient $\bm Q(x)$.

\medskip

As an illustration, in Figure \ref{OU3} we display a sample path of one example of this kind of processes. The example is conveniently chosen so that it points out the effect of the diffusion coefficients.

We take $N=3$ phases and $S=\mathbb{R}$. The simulation of these processes are made according to the stochastic differential equation
$$
dX_t=\tau_{Y_t}(X_t)dt+\sigma_{Y_t}(X_t)dB_t,\quad Y_t=1,2,3,
$$
where
$$
\sigma_i^2(x)=i^2,\quad\tau_i(x)=-ix,\quad i=1,2,3.
$$
Here $B_t$ denotes the standard Brownian motion.

The three phases evolve as Orstein-Ulenbeck processes with different parameters. The diffusion coefficients grow with the number of phase and the drift coefficients are more intense with the number of phase. The transition of phases goes according to the matrix
$$
\bm Q(x)=\begin{pmatrix}
           -\frac{2+x^2}{1+x^2} & \frac{5+3x^2}{6(1+x^2)} & \frac{7+3x^2}{6(1+x^2)} \\
          1+\frac{x^2}{400} & -2-\frac{x^2}{100} & 1+\frac{3x^2}{400} \\
           \frac{1}{2}+\frac{1}{2}e^{-x^2} & \frac{1}{2}+\frac{1}{2}e^{-x^2} & -1-e^{-x^2} \\
         \end{pmatrix}.
$$
In Figure \ref{OU3} we clearly see the evolution of the particle at each phase, starting at $x=0$ and phase 2. We observe that the intensity of the process is higher with the number of phase.

%The transition of phases and the waiting times at each phase are quite homogeneous, except when the particle is far from 0, where it tends to spend more time in phase 2.

%\begin{figure}[h]
%\begin{center}
%\vspace{-7.0cm}
%\includegraphics[height=20cm]{Hermite1.pdf}
%\vspace{-8.0cm}
%\end{center}
%\caption{Sample path of the example introduced in Section XX for the special value of $C=1/2$, starting at phase 2, with four changes, marked with %vertical lines. We observe that oscillations in phase 2 are much closer than oscillations in phase 1.}
%\label{Hermite1}
%\end{figure}

%\begin{figure}[h]
%\begin{center}
%\vspace{-7.0cm}
%\includegraphics[height=20cm]{OneStepN=3.pdf}
%\vspace{-8.0cm}
%\end{center}
%\caption{Sample path of the example introduced in Section XX for the special values of $\alpha=0, \beta=1, k=1.9$, starting at phase 1, with five %changes, marked with vertical lines. We observe clearly the effect of each of the phases.}
%\label{OneStep}
%\end{figure}

\section{Differential and spectral properties of $\mathcal{A}$}\label{NEW}

In this section we will show the main theoretical results of this paper. We will start deriving the corresponding matrix-valued backward and forward equations associated with these processes. Then we will give an explicit formula of the matrix-valued transition density $\bm{P}(t;x,y)$ defined in (\ref{TPD}) in terms of the matrix-valued eigenfunctions of the infinitesimal generator $\mathcal{A}$ defined in (\ref{ForIO}). Finally, we will study some differential equations associated with some functionals related with the recurrence and the invariant distribution of these processes.

\subsection{The matrix-valued backward and forward equations}

For all $\bm{f}(x)\in\mathfrak{B}(S^{N})$ that are twice continuous differentiable we have that the matrix-valued expectation defined in \eqref{MVE} satisfies the backward differential equation
\begin{equation}\label{BDE}
 \frac{\partial}{\partial t}\bm{T}_t\bm{f}(x)=\mathcal{A}\bm{T}_t\bm{f}(x),
\end{equation}
with initial condition $\bm{T}_0\bm{f}(x)=\bm{f}(x)$, where $\mathcal{A}$ is defined in \eqref{ForIO}. This is easy to see by definition of $\mathcal{A}$ in \eqref{infgen} applied to $\bm{f}=\bm{T}_t\bm{f}$, using the Chapman-Kolmogorov equation \eqref{CKEQ}. The initial condition is a consequence of the assumption \eqref{AS3}. We also have that $\bm{T}_t$ and $\mathcal{A}$ commute on all functions in the domain of $\mathcal{A}$.

As a consequence, the transition density $\bm P(t;x,y)$ satisfies the \emph{Kolmogorov backward equation}
\begin{equation}\label{BDETD}
  \frac{\partial \bm P(t;x,y)}{\partial t} =\frac{1}{2}\bm{A}(x)\frac{\partial^2\bm P(t;x,y)}{\partial x^2}+\bm{B}(x)\frac{\partial\bm P(t;x,y)}{\partial x}+\bm{Q}(x)\bm P(t;x,y),
\end{equation}
applicable for $t>0$ and $x,y\in S$.

\bigskip

The  \emph{matrix-valued forward, evolution} or \emph{Fokker-Planck equation} for $\bm P(t;x,y)$ can be regarded as dual to \eqref{BDETD} and the pertinent variables are $t$ and $y$ (the state variable at time $t$ rather than the initial value $x$).

For any integrable vector-valued function $\bm g$, define the operator
\begin{equation*}\label{Tradj}
    \bm{T}_t^* \bm{g}(y)=\int_S\bm{P}^*(t;x,y)\bm{g}(x)dx.
\end{equation*}
Observe that the integration of this operator is with respect to $x$, rather than $y$, as in \eqref{MVE}. This operator is \emph{adjoint} to $\bm{T}_t$ in the sense that
\begin{align}\label{adjT}
\langle\bm{T}_t\bm f,\bm g\rangle&=\int_S \bm g^*(x)\left(\int_S\bm P(t;x,y)\bm f(y)dy\right)dx=\int_S\int_S \bm g^*(x)\bm P(t;x,y)\bm f(y)dxdy\\
\nonumber &=\int_S \left(\int_S\bm P^*(t;x,y)\bm g(x)dx\right)^*\bm f(y)dy=\langle\bm f,\bm{T}_t^*\bm g\rangle.
\end{align}
Here $\langle\bm f,\bm g\rangle=\int_S \bm g^*(x)\bm f(x)dx$  is a real-valued inner product for any two integrable vector-valued functions. Observe that in case we are using matrix-valued functions (as we will do in Section \ref{SRSEC}) we have to use the same inner product $\langle\bm F,\bm G\rangle=\int_S \bm G^*(x)\bm F(x)dx$, but now this inner product is matrix-valued.

If $\bm f$ is twice continuously differentiable and vanishes outside a compact subset of $S$, then one may differentiate with respect to $t$ in \eqref{adjT} and interchange the orders of integration and differentiation to get, using \eqref{BDE},
\begin{equation}\label{AdjPt}
  \left\langle\bm f,\frac{\partial}{\partial t}\bm{T}_t^*\bm g\right\rangle=\left\langle\frac{\partial}{\partial t}\bm{T}_t\bm f,\bm g\right\rangle=\langle\bm{T}_t\mathcal{A}\bm{f},\bm g\rangle=\langle\mathcal{A}\bm{f},\bm{T}_t^*\bm g\rangle=\langle\bm{f},\mathcal{A}^*\bm{T}_t^*\bm g\rangle.
\end{equation}
Here the formal adjoint $\mathcal{A}^*$ of $\mathcal{A}$ with respect to the inner product $\langle\cdot,\cdot\rangle$ is given by
\begin{equation}\label{AdjA}
(\mathcal{A}^*\bm g)(y)=\frac{1}{2}\frac{\partial^2}{\partial y^2}(\bm A^*(y)\bm g(y))-\frac{\partial}{\partial y}(\bm B^*(y)\bm g(y))+\bm{Q}^*(y)\bm g(y).
\end{equation}
From \eqref{AdjPt} we get the forward differential equation
\begin{equation*}\label{FDE}
 \frac{\partial}{\partial t}\bm{T}_t^*\bm{g}(y)=\mathcal{A}^*\bm{T}_t^*\bm{g}(y),
\end{equation*}
or, equivalently, for the transition density $\bm P(t;x,y)$
\begin{equation*}
  \frac{\partial \bm P^*(t;x,y)}{\partial t} =\frac{1}{2}\frac{\partial^2}{\partial y^2}(\bm A^*(y)\bm P^*(t;x,y))-\frac{\partial}{\partial y}(\bm B^*(y)\bm P^*(t;x,y))+\bm{Q}^*(y)\bm P^*(t;x,y).
\end{equation*}
In other words, taking adjoint
\begin{equation}\label{FDETD}
  \frac{\partial \bm P(t;x,y)}{\partial t} =\frac{1}{2}\frac{\partial^2}{\partial y^2}(\bm P(t;x,y)\bm A(y))-\frac{\partial}{\partial y}(\bm P(t;x,y)\bm B(y))+\bm P(t;x,y)\bm{Q}(y).
\end{equation}

Observe that the forward differential equation \eqref{FDETD} has the coefficients of the infinitesimal operator multiplied \emph{on the right}, while the backward differential equation \eqref{BDETD} \emph{on the left}.

%From the Chapman-Kolmogorov equation \eqref{CKEQ}, i.e.
%\begin{equation*}
% \bm{P}(t+s;x,y)=\int_S\bm{P}(t;x,z)\bm{P}(s;z,y)dz,
%\end{equation*}
%and differentiating both sides with respect to $s$ and using the backward equation \eqref{BDETD} satisfied by $\bm{P}(s;z,y)$ we obtain (as a consequence of the semigroup property
%\begin{align*}
%\frac{\partial\bm{P}(t+s;x,y)}{\partial t}&=\frac{\partial\bm{P}(t+s;x,y)}{\partial s}\\ &=\int_S\bm{P}(t;x,z)\left[\frac{1}{2}\bm{A}(z)\frac{\partial^2\bm{P}(s;z,y)}{\partial z^2}+\bm{B}(z)\frac{\partial\bm{P}(s;z,y)}{\partial z}+\bm{Q}(z)\bm{P}(s;z,y)\right]dz,
%\end{align*}
%Now, integration by parts, and assuming the contribution from the boundaries vanish, transforms the formula above into
%\begin{equation*}
% \frac{\partial\bm{P}(t+s;x,y)}{\partial t}=\int_S\left\{\frac{\partial^2}{\partial z^2}\left[\frac{1}{2}\bm{P}(t;x,z)\bm{A}(z)\right]-\frac{\partial}{\partial z}[\bm{P}(t;x,z)\bm{B}(z)]+\bm{P}(t;x,z)\bm{Q}(z)\right\}\bm{P}(s;z,y)dz.
%\end{equation*}
%Since $\bm{P}(s;z,y)$ approaches the delta measure concentrating at $y$, then, as $s\rightarrow0$ the above relation passes into
%\begin{equation*}
% \frac{\partial\bm{P}(t;x,y)}{\partial t}=\frac{\partial^2}{\partial y^2}\left[\frac{1}{2}\bm{P}(t;x,y)\bm{A}(y)\right]-\frac{\partial}{\partial y}[\bm{P}(t;x,y)\bm{B}(y)]+\bm{P}(t;x,y)\bm{Q}(y),
%\end{equation*}
%which is the \emph{forward equation}.

\bigskip

In principle we can try to solve explicitly the backward (or forward) differential equation with specific boundary values to obtain a general expression of the transition probability density $\bm P(t;x,y)$. Nevertheless, this task is not straightforward at all.

There is a very simple case, when $\bm A(x)=\sigma^2(x)\bm I$, $\bm B(x)=\tau(x)\bm I$ and $\bm Q(x)=\bm Q$, i.e., a constant matrix. In this case it is easy to see that
$$
\bm P(t;x,y)=p(t;x,y)e^{t\bm Q},
$$
where $p(t;x,y)$ is the transition probability density of the one dimensional diffusion process with drift $\tau(x)$ and diffusion coefficient $\sigma^2(x)$. In this case both the continuous and discrete components are \emph{independent} and the transition of phases and the lifetime spent on them is independent of the position of $X_t$.

We can try to generalize this separation of variables in the same case as before, but with $\bm Q(x)$ depending on $x$. In this case substituting $\bm P(t;x,y)=p(t;x,y)e^{t\bm Q(x)}$ in \eqref{BDETD} leads to a nonlinear second-order differential equation satisfied by $\bm Q(x)$, which is not straightforward to solve. Things get more complicated assuming $\bm A(x)$ and $\bm B(x)$ diagonal with different entries. Even when $\bm Q(x)$ is constant it is not easy to solve \eqref{BDETD}.

%We can use other approaches like Fourier analysis, to analyze the difficulties of the problem. For instance, let us try to solve the forward equation for the simple case when $\bm A(x)=\bm A$ is a constant diagonal matrix, $\bm B(x)=\bm 0$ and $\bm Q(x)=\bm Q$ is also constant. The forward equation is then
%$$
%  \frac{\partial \bm P(t;x,y)}{\partial t} =\frac{1}{2}\frac{\partial^2}{\partial y^2}\bm P(t;x,y)\bm A+\bm P(t;x,y)\bm{Q}.
%$$
%Denoting $\hat{\bm P}(t;x,z)$ the Fourier transform of $\bm P(t;x,y)$ (entry by entry) leads to (see pp. 390 of \cite{BW})
%$$
%\hat{\bm P}(t;x,z)=e^{ixz}e^{-\frac{1}{2}tz^2\bm A+t\bm Q}.
%$$
%Since $\bm A$ and $\bm Q$ do not commute in general the right hand part of the formula above is not straightforward and we have problems in calculating the inverse of the Fourier transform.

\medskip

Therefore, we need different approaches of obtaining or approximating the transition probability density $\bm P(t;x,y)$. One of them, as we will see in the next section, is via the eigenfunctions of the infinitesimal operator $\mathcal{A}$.

\subsection{The spectral representation of the matrix-valued transition density}\label{SRSEC}

As we have seen in the last section, the explicit computation of $\bm P(t;x,y)$ can be something difficult to handle. In this section we show one way of obtaining a formula for $\bm P(t;x,y)$ in terms of the matrix-valued eigenfunctions of the infinitesimal operator $\mathcal{A}$.

Consider now the space $\LWd$ of all matrix-valued functions $\bm{F}(x)$ such that
\begin{equation*}\label{innerW}
    \langle \bm F,\bm F\rangle_{\bm W}=\int_S\bm F^*(x)\bm W(x)\bm F(x)dx<\infty,
\end{equation*}
where $d\bm W$ is a \emph{weight matrix} with a smooth density $\bm W$ with respect to the Lebesgue measure (see \cite{DG1} for a more concise definition of a weight matrix). In the above definition we mean that the integral is finite entry by entry. This induces a matrix-valued inner product for any two matrix-valued functions $\bm F,\bm G\in\LWd$, denoted by
\begin{equation}\label{inner}
  \langle \bm F,\bm G\rangle_{\bm W}= \int_S\bm G^*(x)\bm W(x)\bm F(x)dx.
\end{equation}
This is not an inner product in the common sense, but it has properties similar to the usual scalar inner products. It is also possible to define a scalar product between two matrix-valued functions (see \cite{DPS}), given by
$$
(\bm F,\bm G)=\mbox{Tr}\left(\langle \bm F,\bm G\rangle_{\bm W}\right).
$$
Therefore, $\LWd$ with the norm $\|\bm F\|_{\bm W}=\mbox{Tr}\left(\langle \bm F,\bm F\rangle_{\bm W}\right)^{1/2}$ is a Hilbert space and (\ref{inner}) is the inner product (in fact, the set of equivalence classes $\bm F\sim \bm G$ if $\|\bm F-\bm G\|_{\bm W}=0$).

\bigskip

The idea behind this method is to find a representation of the transition operator $\bm{T}_t\bm{F}(x)$ (see \eqref{MVE}) using the method of separation of variables as a superposition (linear combination) involving eigenfunctions of the infinitesimal operator $\mathcal{A}$.

Suppose that we know a set of countable normalized matrix-valued eigenfunctions $(\bm\Phi_n(x))_n$ of $\mathcal{A}$ with corresponding eigenvalues $\bm\Gamma_n$ (Hermitian), i.e. $\mathcal{A}\bm\Phi_n(x)=\bm\Phi_n(x)\bm\Gamma_n$. Normalization means $\langle \bm\Phi_n,\bm\Phi_m\rangle_{\bm W}=\delta_{nm}\bm I$ for certain weight matrix $\bm W$.

If the set of finite linear combinations of eigenfunctions is complete in $\LWd$, then each $\bm F\in\LWd$ has a Fourier expansion of the form
\begin{equation*}
 \bm F(x)=\sum_{n=0}^{\infty}\bm\Phi_n(x)\langle \bm F,\bm\Phi_n\rangle_{\bm W}.
\end{equation*}
Consider the superposition defined by
\begin{equation*}
 \bm U_{\bm F}(t,x)=\sum_{n=0}^{\infty}\bm\Phi_n(x)e^{t\bm\Gamma_n}\langle \bm F,\bm\Phi_n\rangle_{\bm W}.
\end{equation*}
Here the exponential for matrix-valued argument is defined using the usual Taylor expansion. It is easy to see that $\bm U_{\bm F}(t,x)$ satisfies the backward differential equation
\begin{equation*}\label{BDES}
 \frac{\partial}{\partial t}\bm U_{\bm F}(t,x)=\mathcal{A}\bm U_{\bm F}(t,x),
\end{equation*}
with initial condition $\bm U_{\bm F}(0,x)=\bm{F}(x)$. So if there is uniqueness for a sufficiently large class of initial functions $\bm F$ then we get
$$
\bm T_t\bm F(x)=\bm U_{\bm F}(t,x),
$$
which means
$$
\int_S\bm P(t;x,y)\bm F(y)dy=\int_S\left(\sum_{n=0}^{\infty}\bm\Phi_n(x)e^{t\bm\Gamma_n}\bm\Phi_n^*(y)\bm W(y)\right)\bm F(y)dy.
$$
Therefore, $\bm P(t;x,y)$ is given by
\begin{equation}\label{SRPD}
 \bm P(t;x,y)=\sum_{n=0}^{\infty}\bm\Phi_n(x)e^{t\bm\Gamma_n}\bm\Phi_n^*(y)\bm W(y).
\end{equation}

\begin{nota}
It is well known that $\bm\Psi_n(x)=\bm\Phi_n(x)\bm U_n$ with $\bm U_n^*\bm U_n=\bm I$ is always another family of normalized eigenfunctions with eigenvalues $\bm U_n^*\bm\Gamma_n\bm U_n$. However, it is straightforward to see that the spectral representation \eqref{SRPD} of  $\bm P(t;x,y)$ is invariant under this change.
\end{nota}

\bigskip

Observe that the condition of having eigenfunctions with Hermitian eigenvalues forces the infinitesimal operator $\mathcal{A}$ to be \emph{self-adjoint} with respect to the inner product \eqref{inner}, i.e. $\langle \mathcal{A}\bm F,\bm G\rangle_{\bm W}=\langle \bm F,\mathcal{A}\bm G\rangle_{\bm W}$ for all $\bm F, \bm G\in\LWd$. This means, using the formal definition of the adjoint of $\mathcal{A}$ given in \eqref{AdjA}, that the coefficients of the infinitesimal operator $\mathcal{A}$ are subject to certain \emph{symmetry equations} with respect to the weight matrix $\bm W$ (with the corresponding boundary conditions)
\begin{equation}\label{symeqs}
\begin{split}
\bm A^*(x)\bm W(x)&=\bm W(x)\bm A(x),\\
\bm B^*(x)\bm W(x)&=(\bm W(x)\bm A(x))'-\bm W(x)\bm B(x),\\
\bm Q^*(x)\bm W(x)&=\frac{1}{2}(\bm W(x)\bm A(x))''-(\bm W(x)\bm B(x))'+\bm W(x)\bm Q(x).
\end{split}
\end{equation}
These symmetry equations were already derived in \cite{DG1, GPT5}. The goal in those papers were to find matrix-valued orthogonal polynomials satisfying second-order differential equations of Sturm-Liouville type (i.e. the differential coefficients are matrix polynomials of degree less than or equal to the order of differentiation). Observe that our infinitesimal operator $\mathcal{A}$ does not have to be in principle of Sturm-Liouville type. In fact, we are more interested in examples with $\bm Q(x)$ depending on $x$, along with diagonal coefficients $\bm A(x)$ and $\bm B(x)$. In this line, we are closer to differential operators as the ones introduced in \cite{GPT1}. Later it was shown (see \cite{GPT6}) that the operators introduced in \cite{GPT1} are intimately related to second-order differential equations of Sturm-Liouville type.

\subsection{Differential equations associated with some functionals}

In this section we will study some differential equations associated with some functionals related with the recurrence and the invariant distribution of these processes.

\medskip

The \emph{hitting time} of a one dimensional diffusion process $\{X_t, 0\leq t<\zeta\}$ to the state $z$ is defined by $T_z=\inf\{t\geq0: X_t=z\}$ or $T_z=\infty$ if $X_t\neq z$ for $0\leq t<\zeta$. We use the notation
$$
T^*=\min\{T_c,T_d\},
$$
for the hitting time to $c$ or $d$, the first time $X_t=c$ or $X_t=d$, where $c,d\in S=(a,b)$. For processes starting at $X_0=x\in(c,d)$, this is the same as the \emph{exit time} of the interval $(c,d)$. We have to accommodate this definition with our bivariate Markov processes. The solution will be introducing matrix-valued analysis.

We will assume that our process is \emph{regular} in the interior of $S$, i.e.
$$
\mbox{Pr}\{T_y<\infty, Y_{T_y}=j|X_0=x, Y_0=i\}>0,
$$
for any $x,y$ in the interior of $S$ and any phases $i,j$.

Consider the matrix-valued function $\bm{U}(x)=(U_{ij}(x))$, where every entry $(i,j)$ is given by
\begin{equation}\label{UU}
U_{ij}(x)=\mbox{Pr}\{T_d<T_c, Y_{T_d}=j|X_0=x, Y_0=i\},\quad c<x<d,
\end{equation}
i.e., starting at $x\in(c,d)$ and phase $i$, $U_{ij}(x)$ is the probability that the process reaches $d$ before $c$ and at that time the process is in phase $j$.

\begin{teor}\label{THU}
 The matrix-valued function $\bm{U}(x)$ defined in (\ref{UU}) possesses two bounded derivatives for $c<x<d$ and satisfies the following differential equation
\begin{equation}\label{DEUU}
 \frac{1}{2}\bm A(x) \bm U''(x)+\bm B(x)\bm U'(x)+\bm Q(x)\bm U(x)=\bm 0,\quad c<x<d,\quad \bm U(c)=\bm 0,\quad \bm U(d)=\bm I.
\end{equation}
\end{teor}
\proof
The proof is essentially the same as in the scalar case. For a heuristic justification see pp. 193 of \cite{KT}. For more technical details the reader should consult Proposition 9.1 of Chapter V in \cite{BW}.
%The boundary conditions are obvious since $U_{ij}(x)$ gives the probability that the process, starting at $x$ and phase $i$, reaches $b$ before $a$ and at that time the process is in phase $j$. Now choose, a time duration $h$ and define $\Delta x=X_h-x$. The matrix-valued version of the law of total probabilities gives
%\begin{equation*}
% \bm{U}(x)=\bm{E}[\bm{U}(X_h)|X_0=x]+\bm{o}(h)=\int_S\bm{P}(t;x,y)\bm{U}(x+\Delta x)dy+\bm{o}(h),
%\end{equation*}
%according to the definition \eqref{MVE}, where the error $\bm{o}(h)$ is a matrix of smaller order than $h$. Expanding in Taylor series $\bm{U}(x+\Delta x)$ of the form
%$$
%\bm{U}(x+\Delta x)=\bm{U}(x)+\Delta x\bm{U}'(x)+\frac{1}{2}(\Delta x)^2\bm{U}''(x)+\cdots
%$$
%and using assumptions \eqref{AS21}, \eqref{AS22} and \eqref{AS3}, we obtain
%$$
%\bm{U}(x)=(\bm I+h\bm{Q}(x))\bm{U}(x)+h\bm{B}(x)\bm{U}'(x)+\frac{1}{2}h\bm{A}(x)\bm{U}''(x)+\bm{o}(h).
%$$
%Subtracting $\bm{U}(x)$ from both sides, dividing by $h$, and letting $h$ decrease to zero we obtain \eqref{DEUU}.
\hfill \endproof
\begin{nota}
In \cite{BW, KT} we can find an explicit solution of this problem in the scalar situation (where $\bm{Q}(x)=\bm 0$) in terms of what they call \emph{scale function} and \emph{speed density} of the process. The presence of $\bm{Q}(x)$, not zero in general, leads to a more complicate system of second-order differential equations with initial conditions
\begin{equation}\label{SDeq}
 \frac{1}{2}\sigma_{i}^2(x)U_{ij}''(x)+\tau_{i}(x)U_{ij}'(x)+Q_{ii}(x)U_{ij}(x)+\sum_{k\neq i}^NQ_{ik}(x)U_{kj}(x)=0, \quad i,j=1,\ldots,N.
\end{equation}
\end{nota}

With the above considerations we can define \emph{recurrence} and \emph{transience} of bivariate Markov processes (in the same line as Definition 9.1 of Chapter V in \cite{BW}). Write the matrix-valued function $\bm R_{x,y}=(R_{x,y})_{i,j}$ where every entry $(i,j)$ is given by
\begin{equation*}\label{UUr}
(R_{x,y})_{i,j}=\mbox{Pr}\{X_t=y, Y_{t}=j|X_0=x, Y_0=i\},\quad x,y\in S,
\end{equation*}
i.e., the probability that starting at $x$ and phase $i$, the process ever reaches $y$ and phase $j$.

For $y>x$, $\bm R_{x,y}$ will be the solution of the differential equation \eqref{DEUU} letting $c\rightarrow a$ and $d\rightarrow y$ and, for $y<x$, $\bm R_{x,y}$ will be the solution of the differential equation \eqref{DEUU} letting $c\rightarrow b$ and $d\rightarrow y$.

\begin{defi}\label{defirec}
A state $y$ is called recurrent if $\bm R_{x,y}\bm e_N=\bm e_N$ for all $x\in S$, where $\bm e_N$ is $N$-dimensional column vector with all entries equal to 1, and transient otherwise. If all states in $S$ are recurrent, then the bivariate Markov process is said to be recurrent.
\end{defi}

\bigskip
Consider now, for a continuous matrix-valued function $\bm G(x)$, finding $\bm{V}(x)$ given by
\begin{equation}\label{VV}
\bm V(x)=\bm E\left[\int_0^{T^*}\bm G(X_s)ds|X_0=x\right],\quad c<x<d,
\end{equation}
where $\bm E$ is defined in (\ref{MVE}).

\begin{teor}\label{THV}
 The matrix-valued function $\bm{V}(x)$ defined in (\ref{VV}) possesses two bounded derivatives for $c<x<d$ and satisfies the following differential equation
\begin{equation}\label{DEVV}
 \frac{1}{2}\bm A(x) \bm V''(x)+\bm B(x)\bm V'(x)+\bm Q(x)\bm V(x)+\bm G(x)=\bm 0,\quad c<x<d,\quad \bm V(c)=\bm V(d)=\bm 0.
\end{equation}
\end{teor}
\proof
As before, the proof is essentially the same as in the scalar case. For a heuristic justification see pp. 194 of \cite{KT}. For more technical details the reader should consult Proposition 10.1 of Chapter V in \cite{BW}.
%The initial conditions are again obvious. Choose a short time $h$. At time $h$ the expectation of the total integral is the expectation of the contribution up to time $h$, which is $(\bm I+h\bm{Q}(x))h\bm{G}(x)+\bm{o}(h)$ since the sample paths and $\bm G(x)$ are continous, plus the expectation of the contribution over the remaining time, which is $\bm{E}[\bm{V}(X_h)|X_0=x]$ using the Markov property. Therefore,
%$$
%\bm{V}(x)=h\bm{G}(x)+\bm{E}[\bm{V}(X_h)|X_0=x]+\bm{o}(h).
%$$
%As before, using Taylor expansions for the second part of the previous right hand formula gives
%$$
%\bm{V}(x)=(\bm I+h\bm{Q}(x))h\bm{G}(x)+(\bm I+h\bm{Q}(x))\bm{V}(x)+h\bm{B}(x)\bm{V}'(x)+\frac{1}{2}h\bm{A}(x)\bm{V}''(x)+\bm{o}(h),
%$$
%and subtracting $\bm{V}(x)$ from both sides, dividing by $h$ and sending $h$ to zero will produce \eqref{DEVV}.
\hfill \endproof

In particular, in the case where
\begin{equation}\label{GGG}
\bm G(x)=\bm e_N\bm e_N^*=\begin{pmatrix}
  1 & 1 & \cdots&1  \\
  1 & 1 & \cdots&1 \\
\vdots&\vdots&\ddots&\vdots\\
1 & 1 & \cdots&1 \\
\end{pmatrix},
\end{equation}
we have that the matrix-valued function $\bm{V}(x)=(V_{ij}(x))$ gives
\begin{equation}\label{Vrec}
 V_{ij}(x)=E\left[T^*, Y_{T^*}=j | X_0=x, Y_0=i\right],
\end{equation}
i.e., starting at $X_0=x$ and phase $i$, $V_{ij}(x)$ is the mean time to reach either $c$ or $d$ and at that time the process is at phase $j$.
\begin{nota}
Again, in \cite{BW, KT} we can find an explicit solution of this problem in the scalar situation in terms of the Green function of the process, but the presence of $\bm{Q}(x)$ will lead to a system of differential equations similar to \eqref{SDeq}, which is more difficult to solve.
\end{nota}

\begin{defi}\label{defiposrec}
Assume our process is recurrent. By letting $c\rightarrow a$ and $d\rightarrow b$ in \eqref{DEVV}, we say that the bivariate Markov process is positive recurrent if all entries of $\bm V(x)$ in \eqref{Vrec}, with $\bm G(x)$ given by \eqref{GGG}, are finite. Otherwise is null recurrent.
\end{defi}

\bigskip

Finally, let us discuss the \emph{invariant distribution} of a bivariate Markov process, which is the behavior of $\bm P(t;x,y)$ as $t\rightarrow\infty$. This distribution should be independent of the initial state $X_0=x$ and the initial phase $Y_0=i$. Therefore we should expect a \emph{row vector-valued} invariant distribution
$$
\bm \psi(y)=(\psi_{1}(y),\psi_{2}(y),\ldots,\psi_{N}(y)),
$$
satisfying
\begin{equation*}\label{InvD}
    \bm \psi(y)=\int_S\bm \psi(x)\bm P(t;x,y)dx,\quad t>0.
\end{equation*}
It is also expected, as in the scalar situation (see Theorem 12.2 of Chapter V in \cite{BW} or pp. 241 of \cite{KT}), that the existence and unicity of such invariant distribution is restricted to the case when the process is positive recurrent.

It is easy to see, mimicking the derivation of the forward differential equation \eqref{FDETD}, that $\bm \psi(y)$ satisfies
\begin{equation}\label{DEID}
    \frac{1}{2}(\bm \psi(y)\bm A(y))''-(\bm \psi(y)\bm B(y))'+\bm \psi(y)\bm Q(y)=\bm 0.
\end{equation}
In order to have a vector-valued probability distribution we must impose that $0\leq\psi_{j}(y)\leq 1$ and
\begin{equation}\label{NormID}
    \left(\int_S\bm \psi(y)dy\right)\bm e_N=1.
\end{equation}

%In the same way as in the scalar case (see pp. 220 of \cite{KT}) we should expect that the vector-valued invariant distribution is approached to the extent that
%$$
%\lim_{t\rightarrow\infty} P_{i,\cdot}(t;x,y)=\bm\psi(y),
%$$
%holds in some appropriate sense for any row $i=1,\ldots,N$. Here $P_{i,\cdot}(t;x,y)$ denotes any row of $\bm P(t;x,y)$.

%If the infinitesimal operator $\mathcal{A}$ is self-adjoint with respect to certain weight matrix $\bm W$ then we will have a spectral representation \eqref{SRPD} of $\bm P(t;x,y)$. Therefore, by looking at the eigenvalues $\bm\Gamma_n$ we can generate an invariant distribution of the process.

If the infinitesimal operator $\mathcal{A}$ is self-adjoint with respect to certain weight matrix $\bm W$ then we can solve directly \eqref{DEID}. This equation resembles the third symmetry equation in \eqref{symeqs}. Substituting
$$
\bm\psi(y)=c\bm e^*_N\bm W(y),\quad c>0,
$$
in \eqref{DEID} and using the third equation in \eqref{symeqs} give
$$
c\bm e^*_N\left[\frac{1}{2}(\bm W(y)\bm A(y))''-(\bm W(y)\bm B(y))'+\bm W(y)\bm Q(y)\right]=c\bm e^*_N\bm Q^*(y)\bm W(y)=\bm 0,
$$
as a consequence of $\bm Q(y)e_N=\bm 0$, for all $y$. The normalization constant $c$ is given by \eqref{NormID}:
$$
c=\left(\int_S\bm e^*_N\bm W(y)\bm e_Ndy\right)^{-1}.
$$
Therefore an explicit formula for the invariant distribution in this case is given by
\begin{equation}\label{IDW}
    \bm\psi(y)=\left(\int_S\bm e^*_N\bm W(y)\bm e_Ndy\right)^{-1}\bm e^*_N\bm W(y).
\end{equation}

\bigskip

\section{A variant of the Wright-Fisher diffusion model involving only mutation effects}\label{EXAM}

The Wright-Fisher diffusion model involving only mutation effects considers a big population of constant size composed of two types A and B. The intensity of mutation from A to B and from B to A are given by two positive constants, $\frac{1+\beta}{2}$ and $\frac{1+\alpha}{2}$, respectively ($\alpha,\beta>-1$). As the size of the population tends to infinity, it is very well known (see pp. 177 of \cite{KT}) that this model can be described by a diffusion process whose state space is the unit interval $S=[0,1]$ with drift and diffusion coefficient given by
\begin{equation}\label{JacCoef2}
    \tau(x)=\alpha+1-x(\alpha+\beta+2),\quad\sigma^2(x)=2x(1-x),\quad\alpha,\beta>-1.
\end{equation}
The variable $x$ is the rate of the number of A-types with respect to the total population. These coefficients correspond with the Jacobi diffusion model, and the probability density can be described in terms of the Jacobi polynomials on $[0,1]$.

In this section we will study a variation of the Wright-Fisher model with only mutation effects given by a bivariate Markov process whose state space is $[0,1]\times\{1,2,\ldots,N\}$, so that we can apply all the analysis introduced in the previous sections. In our case, there will be $N$ possible phases with different drift coefficients (the diffusion coefficient will always be the same). Also, there will be an extra positive parameter $k,$ depending on $\beta$, which measures how the process moves through all the phases. All these considerations will be studied deeply in Sections \ref{PROB} and \ref{FINREM}.

The example we study comes from group representation theory, and it is related with matrix-valued spherical functions associated with the complex projective space. It was first introduced in \cite{GPT1} for special choice of the parameters (see also \cite{GPT6}). Later in \cite{PT1} it was related with matrix-valued orthogonal polynomials. In Section \ref{1step} we will give the coefficients relevant in this example to treat it as a bivariate Markov process. The recurrence coefficients of this example have also been related with certain quasi-birth-and-death processes in \cite{GdI2}. Recently, new applications of these processes have been related to urn or Young diagram models \cite{GPT7}. In both cases the differential operator played no role. In this paper we will give a meaning of the differential operator in terms of bivariate Markov processes.

\subsection{A family of examples arising from the complex projective space}\label{1step}

In what follows we will use $\bm E_{ij}$ to denote the matrix with 1 at entry $(i,j)$ and 0 otherwise. Let $N\in\{1,2,\ldots\}$, $\alpha,\beta>-1$ and $0<k<\beta+1$. $\bm M$ will denote the nilpotent upper triangular matrix
\begin{equation*}\label{AA}
 \bm M=\sum_{i=1}^{N-1}\bm E_{i,i+1},
\end{equation*}
while $\bm J$ and $\breve{\bm J}$ the following diagonal matrices
\begin{equation*}\label{JJ}
 \bm J=\sum_{i=1}^{N}\bm (N-i)\bm E_{ii},\quad  \breve{\bm J}=\sum_{i=1}^{N}\bm (i-1)\bm E_{ii}=(N-1)\bm I-\bm J.
\end{equation*}
Finally, the diagonal matrix $\bm H$ will denote
\begin{equation*}\label{HH}
 \bm H=\sum_{i=1}^{N}\omega_i\bm E_{ii},
\end{equation*}
where
$$
\omega_i=\begin{pmatrix}
                                      \beta-k+i-1 \\
                                       i-1 \\
                                     \end{pmatrix}\begin{pmatrix}
                                       N+k-i-1 \\
                                       N-i \\
                                     \end{pmatrix}.
$$
Observe that the weights $\omega_i$ correspond to the weights of the \emph{Hahn polynomials} supported on $\{1,2,\ldots,N\}$ with parameters $\beta-k>-1$ and $k-1>-1$ (see for instance pp. 345 of \cite{AAR}).

\medskip

The pair $\{\bm W, \mathcal{A}\}$, where
\begin{equation}\label{WW}
 \bm W(x)=x^{\alpha}(1-x)^{\beta}\bm Hx^{\bm J},\quad x\in (0,1),
\end{equation}
and
\begin{equation}\label{DD}
 \mathcal{A}=\frac{1}{2}\bm A(x)\frac{d^2}{dx^2}+\bm B(x)\frac{d}{dx}+\bm Q(x)\frac{d^0}{dx^0},
\end{equation}
where
\begin{align}\label{ABQ}
&\bm A(x)=2x(1-x)\bm I,\quad \bm B(x)=(\alpha+1)\bm I+\bm J-x((\alpha+\beta+2)\bm I+\bm J),\\
\nonumber& \bm Q(x)=\frac{1}{1-x}\left(\bm J\left((\beta-k+1)\bm I+\breve{\bm J}\right)(\bm M-\bm I)+x\left(\breve{\bm J}(k\bm I+\bm J)(\bm M^*-\bm I)\right)\right),
\end{align}
is a \emph{symmetric pair}, meaning that $\mathcal{A}$ is self-adjoint with respect to $\bm W$ (with respect to the inner product \eqref{inner}). See \cite{GPT5, PT1} for definitions.

The \emph{tridiagonal} matrix $\bm Q(x)$ can be written as
\begin{equation}\label{QQ}
\bm Q(x)=\sum_{i=2}^{N}\mu_i(x)\bm E_{i,i-1}-\sum_{i=1}^{N}(\lambda_i(x)+\mu_i(x))\bm E_{ii}+\sum_{i=1}^{N-1}\lambda_i(x)\bm E_{i,i+1},
\end{equation}
where
\begin{equation}\label{Hahnx}
\begin{split}
  \lambda_i(x) =& \frac{1}{1-x}(N-i)(i+\beta-k), \\
\mu_i(x) =& \frac{x}{1-x}(i-1)(N-i+k).
\end{split}
\end{equation}

We remark that $\lambda_i(x)$ and $\mu_i(x)$ are $x$-dependent variations of the coefficients of the second-order difference equation satisfied by the Hahn polynomials (see for instance pp. 346 of \cite{AAR}).

Observe that $\bm A(x)$ and $\bm B(x)$ are diagonal matrices and the diagonal entries of $\bm A(x)$ are positive for all $x\in (0,1)$. Also $\lambda_i(x)$ and $\mu_i(x)$ are always positive and $\bm Q(x)\bm e_N=\bm 0$ for every $x\in (0,1)$. Therefore $\bm Q(x)$ is the intensity matrix of a classical continuous time Markov chain with finite state space. All these requirements fit perfectly with the definition of bivariate Markov processes introduced in Section \ref{PRE}.

\medskip

We have many formulas and properties associated with this example at our disposal in the literature. In particular we can find a set of normalized eigenfunctions $(\bm\Phi_n(x))_n$ of the differential operator $\mathcal{A}$ \eqref{DD} with the corresponding eigenvalues $(\bm\Gamma_n)_n$. This will allow us to get a formula for the matrix-valued transition density $\bm P(t;x,y)$ of this process using formula \eqref{SRPD}.

The eigenfunctions $(\bm\Phi_n(x))_n$ of $\mathcal{A}$ (also known as the matrix-valued spherical functions associated with the complex projective space, see \cite{GPT1}) are very close related to matrix-valued orthogonal polynomials and they have been studied in several papers. We use the family of matrix-valued orthogonal polynomials given in \cite{GdI2}, since it contains most of the structural formulas we need for this paper. First we describe the relation between the symmetric pair $\{\bm W, \mathcal{A}\}$ given in \eqref{WW} and \eqref{DD} and the symmetric pair $\{ \widetilde{\bm W}, \mathcal{\widetilde{A}}\}$ given in \cite{GdI2}. This transformation was introduced for the first time in \cite{GPT6} (see also \cite{PT1}). Denote $\bm\Psi(x)$ the following matrix
\begin{equation*}\label{Psi}
 \bm\Psi(x)=e^{\breve{\bm J}\bm M^*}(1-x)^{\breve{\bm J}} \bm T,
\end{equation*}
where $\bm T$ is the constant upper triangular matrix given at the beginning of Section 7 in \cite{GdI2}.

The relation between the pair $\{\bm W, \mathcal{A}\}$ and $\{\widetilde{\bm W}, \mathcal{\widetilde{A}}\}$ is the following
\begin{equation}\label{Relat}
 \widetilde{\bm W}(x)=\bm\Psi^*(x)\bm W(x)\bm\Psi(x),\quad \mathcal{\widetilde{A}}\bm F(x)=\bm\Psi^{-1}(x)\mathcal{A}\left(\bm\Psi(x)\bm F(x)\right).
\end{equation}

Take now the family of matrix-valued orthogonal polynomials introduced in \cite{GdI2}, which we will denote by $(\bm R_n(x))_n$ (in \cite{GdI2} they are denoted by $(Q_n^*(x))_n$). Then the family $(\bm R_n(x))_n$ is eigenfunction of the differential operator $\mathcal{\widetilde{A}}$, i.e. $\mathcal{\widetilde{A}}\bm R_n(x)=\bm R_n(x)\bm\Gamma_n$, where $\bm\Gamma_n$ is the \emph{diagonal} matrix
\begin{equation}\label{Lamb}
\bm\Gamma_n=-n^2\bm I-n((\alpha+\beta+N)\bm I+\breve{\bm J})-\breve{\bm J}((\alpha+\beta-k+1)\bm I+\breve{\bm J}).
\end{equation}
This family has very special properties. In particular it is very easy to normalize since the norms $\|\bm R_n\|_{\widetilde{\bm W}}^2$, given in formula (7.3) of \cite{GdI2}, are \emph{diagonal} matrices.

\medskip

Therefore, define the family
\begin{equation*}\label{Phi}
 \bm\Phi_n(x)=\bm\Psi(x)\bm R_n(x)\|\bm R_n\|_{\widetilde{\bm W}}^{-1}.
\end{equation*}
We have that this family is eigenfunction of the infinitesimal generator $\mathcal{A}$, i.e. $$\mathcal{A}\bm \Phi_n(x)=\bm \Phi_n(x)\bm\Gamma_n,$$ with $\bm\Gamma_n$ defined in \eqref{Lamb}, and orthonormal with respect to $\bm W(x)$ \eqref{WW}. Observe that $(\bm\Phi_n(x))_n$ are also matrix-valued polynomials, but not of degree exactly $n$ (more precisely of degree $n+N-1$).

Therefore, according to formula \eqref{SRPD}, we have
\begin{equation}\label{SRPDex1}
 \bm P(t;x,y)=\sum_{n=0}^{\infty}\bm\Phi_n(x)e^{t\bm\Gamma_n}\bm\Phi_n^*(y)\bm W(y),
\end{equation}
or, in terms of the pair $\{\widetilde{\bm W}, \mathcal{\widetilde{A}}\}$,
\begin{equation*}\label{SRPDex12}
 \bm P(t;x,y)=\bm \Psi(x)\left(\sum_{n=0}^{\infty}\bm R_n(x)e^{t\bm\Gamma_n}\|\bm R_n\|_{\widetilde{\bm W}}^{-2}\bm R_n^*(y)\widetilde{\bm W}(y)\right)\bm\Psi^{-1}(y).
\end{equation*}

\subsection{Probabilistic interpretation}\label{PROB}

Consider the Wright-Fisher diffusion model involving only mutation effects for a big population of constant size composed of two types A and B. The intensity of mutation from A to B and from B to A are described at the beginning of Section \ref{EXAM} and are given in terms of the parameters $\alpha,\beta>-1$ of the example introduced in Section \ref{1step}. Now there will be an indicator denoted by the extra parameter $k$, depending on the intensity $\beta$, which helps the population of A to survive against the population of B. We will give more details in Section \ref{FINREM}.

\medskip

First we remark the importance of having a representation of the matrix-valued probability density $\bm P(t;x,y)$ given by \eqref{SRPDex1}, which can be very well approximated by the first eigenfunctions $\bm \Phi_n(x)$ and the corresponding eigenvalues $\bm\Gamma_n$. We can calculate very accurately what is the probability of reaching any Borel set $B\subset (0,1)$ and certain phase $Y_t=j$ in time $t$ given that initially we start from state $X_0=x$ and phase $Y_0=i$. As an illustration, for $N=4$ phases, starting at $x=1/2$, at time $t=1$ we have
$$
\mbox{Pr}\left\{X_1\in\left(\frac{3}{4},1\right)\bigg|X_0=\frac{1}{2}\right\}=\begin{pmatrix}
  0.12410905 & 0.08138740 & 0.08920446 & 0.1633878\\
  0.11006872 & 0.07334748 & 0.08181737 & 0.1527569\\
  0.09668381 & 0.06555764 & 0.07453306 & 0.1420720\\
  0.08394494 & 0.05801744 & 0.06735379 & 0.1313385\\
\end{pmatrix},
$$
where the entries $(i,j)$ represent the initial phase $Y_0=i$ and final phase $Y_1=j$. Here we have taken the values of the parameters $\alpha=\beta=0$ and $k=\frac{1}{2}$. In this case only 3 eigenfunctions are needed in order to get a very accurate approximation of the probabilities.

\medskip

Our family of bivariate Markov processes $(X_t,Y_t)$ with state space $[0,1]\times\{1,2,\ldots,N\}$ can be described in terms of the stochastic differential equation
$$
dX_t=\tau_{Y_t}(X_t)dt+\sigma(X_t)dB_t,
$$
where (see \eqref{ABQ})
\begin{equation}\label{JacCoefBMP}
    \tau_i(x)=\alpha+1+N-i-x(\alpha+\beta+2+N-i),\quad\sigma^2(x)=2x(1-x),\quad Y_t=i=1,\ldots,N.
\end{equation}
The description of how the process moves through the different phases is given by $\bm Q(x)$ \eqref{QQ}. Observe that $\bm Q(x)$ is tridiagonal, meaning that transitions of phases are only allowed between adjacent phases (like a birth-and-death process).

The behavior of $X_t$ depends on the phase $Y_t$ and the transition of phases and the waiting times at each phase $Y_t$ depend also on the position of $X_t$. About the continuous part $X_t$ we are interested in studying the behavior at the boundaries values, in this case 0 and 1. On the contrary, from the discrete part $Y_t$, we will be interested in analyzing how the process moves through the different phases, including the waiting times at each phase. The behavior of the bivariate Markov processes $(X_t,Y_t)$ will be a combination of both aspects.

\medskip

The diffusion coefficient $\sigma^2(x)$ is always the same, but the drift coefficients $\tau_i(x)$ depend on the phase $i$. For a given phase $i$, moving to the next phase $i+1$ means decreasing one unity of the intensity of mutation $B\rightarrow A$ (in terms of the parameter $\alpha$), and moving to the previous phase $i-1$ means increasing one unity of the same intensity of mutation $B\rightarrow A$. The limiting situation is at phase $N$, where we recover the regular Jacobi diffusion \eqref{JacCoef2}. The intensity of mutation $A\rightarrow B$ (in terms of the parameter $\beta$) is not affected in these coefficients.

In page 239 of \cite{KT} one can find a classification of the boundary states 0 and 1 of the Wright-Fisher model involving only mutation pressures for Jacobi diffusions in terms of the parameters. Both states 0 and 1 are \emph{regular} or \emph{absorbing} boundaries when $-1<\alpha,\beta<0$, respectively. While both states 0 and 1 are \emph{entrance} or \emph{reflecting} boundaries when $\alpha,\beta\geq0$.\footnote{A reflecting boundary cannot be reached from the interior of the state space, but it is possible to consider processes that begin there. An absorbing boundary can both be entered and left.}

\begin{figure}[h]
\begin{center}
\vspace{-5.0cm}
\includegraphics[height=20cm]{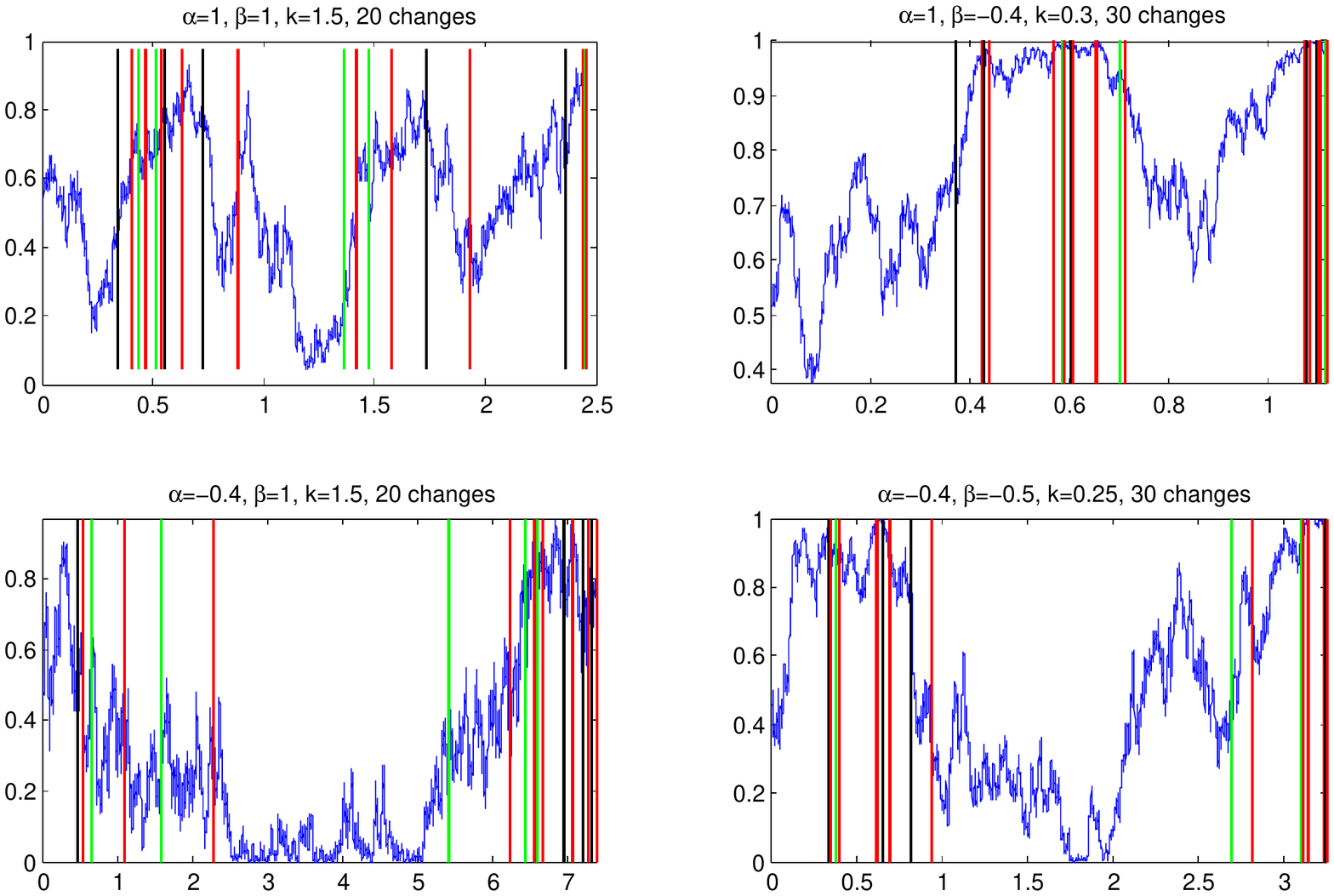}
\vspace{-6.0cm}
\end{center}
\caption{Different behaviors of the boundaries 0 and 1 for different values of $\alpha$, $\beta$ and $k$. In the first plot 0 and 1 are reflecting. In the second plot, 1 is absorbing and 0 reflecting. In the third plot, 1 is reflecting and 0 is absorbing (only at phase 3). Finally in the fourth plot, 0 (only at phase 3) and 1 are absorbing.}
\label{WFB}
\end{figure}

For our bivariate Markov process observe that the boundary state 1 preserves this behavior for any phase. But the boundary state 0 only preserves this behavior at phase $N$. For the rest of phases $i\in\{1,\ldots,N-1\}$ the boundary state 0 is always reflecting since $\alpha+N-i>0$ (see \eqref{JacCoefBMP}). Therefore
\begin{equation*}
1\begin{cases}
\text{is an absorbing boundary } & \text{for } -1<\beta<0, \\
\text{is a reflecting boundary } & \text{for } \qquad\beta\geq0.
\end{cases}
\end{equation*}
\begin{equation*}
0\begin{cases}
\text{is an absorbing boundary } & \text{for } -1<\alpha<0 \quad \text{AND}\quad \text{phase $N$,} \\
\text{is a reflecting boundary } & \text{for }  \qquad\alpha\geq0.
\end{cases}
\end{equation*}
In Figure \ref{WFB} we can check this behavior for four different situations. In all of them we have taken $N=3$ phases. In this figure we have not marked with numbers the transitions of phases (in vertical lines) since we are only interested in the behavior at the boundaries 0 and 1. In the first plot the parameters are taken so that the boundaries 0 and 1 are both reflecting, meaning that the sample path is never going to reach 0 or 1. In the second plot, 1 is absorbing, while 0 is reflecting. We observe that 1 can be reached in principle at any phase. In the third plot, on the contrary, 0 is absorbing while 1 is reflecting. In this case the only phase where the sample path reaches 0 is \emph{only on the last third phase}. In the rest of phases, 0 becomes reflecting. Finally, in the fourth plot, the parameters are taken so that 0 and 1 are both absorbing. The process can reach boundary 1 at any phase, but, again, boundary 0 only when the process is at the third phase.

According to Definition \ref{defirec} our process is always \emph{recurrent} for $\alpha,\beta\geq0$, while \emph{transient} otherwise. This fact has been checked by looking extensively at the numerical solutions of $\bm R_{x,y}$ via the second-order differential equation \eqref{DEUU}. In a similar way,
we can state, using Definition \ref{defiposrec}, that our process is always \emph{positive recurrent} for $\alpha,\beta\geq0$ by looking the numerical solutions of $\bm V(x)$ via the second-order differential equation \eqref{DEVV} with $\bm G(x)$ given by \eqref{GGG}.

\medskip

Let us now study a couple of aspects from the discrete component $Y_t$. First we will study the \emph{waiting times} at each phase depending on the parameters and secondly we will study the \emph{tendency} of moving forward or backward in phases. All these facts depend on the position of $X_t$.

For the waiting times we have to look at the diagonal entries of $\bm Q(x)$ given by (see \eqref{QQ} and \eqref{Hahnx})
\begin{equation*}
\begin{cases}
 Q_{11}(x)=-\frac{1}{1-x}(N-1)(\beta-k+1)& \text{for } i=1, \\
 Q_{ii}(x)=-\frac{1}{1-x}\left[(N-i)(i+\beta-k)+x(i-1)(N-i+k)\right]& \text{for } i=2,\ldots,N-1, \\
 Q_{NN}(x)=-\frac{x}{1-x}(N-1)k& \text{for }  i=N.
\end{cases}
\end{equation*}

\begin{figure}[h]
\begin{center}
\vspace{-5.0cm}
\includegraphics[height=20cm]{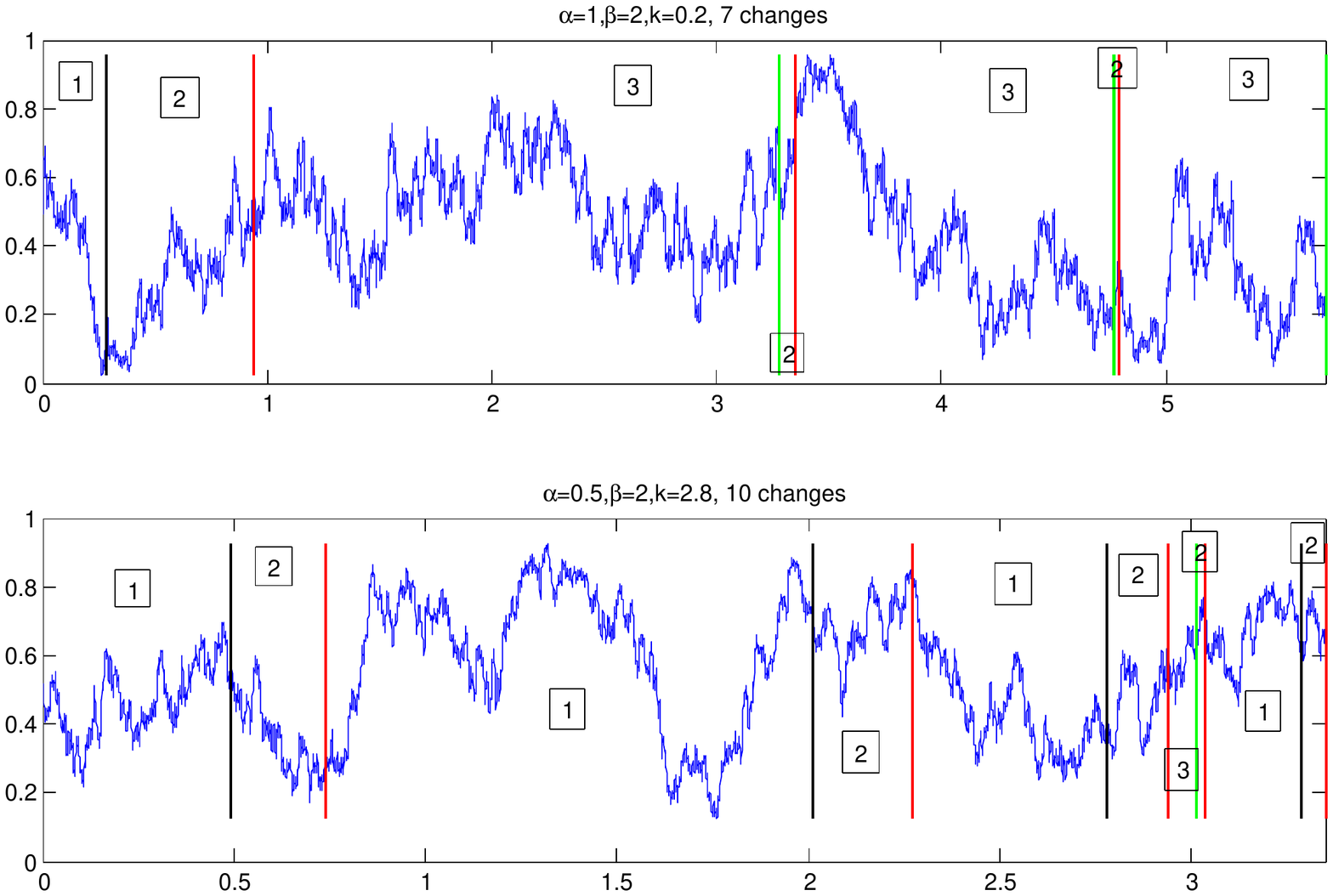}
\vspace{-6.0cm}
\end{center}
\caption{Different behaviors of the phases 1 and $3$ for different values of $k$ ($N=3$). In the first plot $k$ is close to 0, so the process tends to spend more time in the last absorbing phase $3$. In the second plot, $k$ is close to $\beta+1$, so the phase 1 is absorbing.}
\label{WFk}
\end{figure}

We have that if $x$ is close to 1, the diagonal coefficients are very large, meaning that all phases are \emph{instantaneous}, i.e., the waiting times at each phase are very short until $x$ is far from 1. This can be checked in the graphs of Figure \ref{WFB}. If $x$ is close to 0, we have that $Q_{NN}(x)$ is very small, meaning that phase $N$ becomes \emph{absorbing}, i.e., if the process enters this phase and the position of $X_t$ is close to 0, then it tends to spend long periods of time there (see the last two plots of Figure \ref{WFB}). Also phase $N$ becomes absorbing (no matter of the position of $X_t$) if the parameter $k$ is close to 0 (see first plot in Figure \ref{WFk}). The middle phases (from 2 to $N-1$) are never absorbing, since $Q_{ii}(x), i=2,\ldots,N-1,$ are never close to 0. Finally, if $k$ is close to $\beta+1$, we have that phase 1 is absorbing (see second plot in Figure \ref{WFk}). For the rest of cases, the value of these diagonal entries and the position of $X_t$ determine the waiting time in that phase.

\medskip

For the tendency we have to study, given that the process is at phase $Y_t=i$, the probability of moving to the next or to the previous phase. Obviously, when the process is at one of the boundary phases 1 or $N$ it can only go forward or backward, respectively. This probability is given by
\begin{eqnarray*}
  \text{Pr}\{Y_t=i\rightarrow Y_t=i+1\} &=& \frac{\lambda_i(x)}{\lambda_i(x)+\mu_i(x)}, \\
  \text{Pr}\{Y_t=i\rightarrow Y_t=i-1\} &=& \frac{\mu_i(x)}{\lambda_i(x)+\mu_i(x)},
\end{eqnarray*}
where the coefficients $\lambda_i(x)$ and $\mu_i(x)$ are given by \eqref{Hahnx}. Here $x$ is the state $X_t$ of the process at the time $t$ of transition of phases. The forward or backward tendency depends on if $\mu_i(x)<\lambda_i(x)$ or $\mu_i(x)>\lambda_i(x)$, respectively. Define the threshold values
\begin{equation}\label{thre}
    x_0(i)=\frac{(N-i)(i+\beta-k)}{(i-1)(N-i+k)},\quad i=2,\ldots,N-1.
\end{equation}
These values are the $x-$solutions of the algebraic equation $\lambda_i(x)=\mu_i(x)$. Therefore, at the moment of changing from phase $Y_t=i$ to the next or previous phase, we will have a \emph{forward tendency} when $x<x_0(i)$, and a \emph{backward tendency} when $x>x_0(i)$. In order to be relevant, this threshold value $x_0(i)$ must be less than 1. Otherwise, there is always a forward tendency.

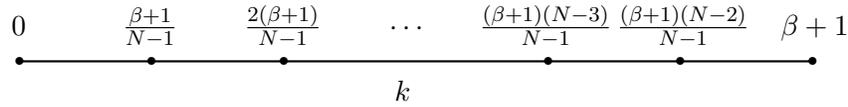
\begin{figure}
\begin{center}
\begin{picture}(200,40)(0,0)
   \linethickness{0.4 mm}
   \thicklines
   \put(-50,0){\line(200,0){300}}
   \put(-50,0){\circle*{3}}
   \put(0,0){\circle*{3}}
   \put(50,0){\circle*{3}}
   \put(150,0){\circle*{3}}
   \put(200,0){\circle*{3}}
   \put(250,0){\circle*{3}}
   \put(-53,10){$0$}
   \put(-11,10){$\frac{\beta+1}{N-1}$}
   \put(35,10){$\frac{2(\beta+1)}{N-1}$}
   \put(90,10){$\cdots$}
   \put(92,-15){$k$}
   \put(124,10){$\frac{(\beta+1)(N-3)}{N-1}$}
   \put(176,10){$\frac{(\beta+1)(N-2)}{N-1}$}
   \put(238,10){$\beta+1$}
   \end{picture}
   \end{center}
   \caption{Range of $k$ for the study or the tendency of the process.}\label{Figkk}
\end{figure}

Let us study the cases for which the threshold values have a meaning. If $x_0(i)>1$ for all $i=1,\ldots,N-1,$ then we will always have a forward tendency. This happens when
$$
k<\frac{(\beta+1)}{N-1}.
$$
Nevertheless, if
\begin{equation}\label{backtend}
k>\frac{(\beta+1)(N-i)}{N-1},
\end{equation}
for some $i$ there can be forward and backward tendency at the same run, depending on the position of $X_t$ at the moment of transition. We have then a maximum backward tendency for $i=2$ in \eqref{backtend}. The range of the values of $k$ depending on the phase $Y_t=i$ are described in Figure \ref{Figkk}. Therefore we have
\begin{equation*}
\begin{cases}
\text{Maximum forward tendency} & \text{for } \quad k<\frac{(\beta+1)}{N-1},\\
\text{Maximum backward tendency}& \text{for } \quad k>\frac{(\beta+1)(N-2)}{N-1},\\
\text{Backward/forward tendency}& \text{otherwise. } \\
\end{cases}
\end{equation*}
This also explains why phases 1 or $N$ are absorbing in terms of the values of $k$.

\medskip

\begin{figure}[h]
\begin{center}
\vspace{-6.0cm}
\includegraphics[height=20cm]{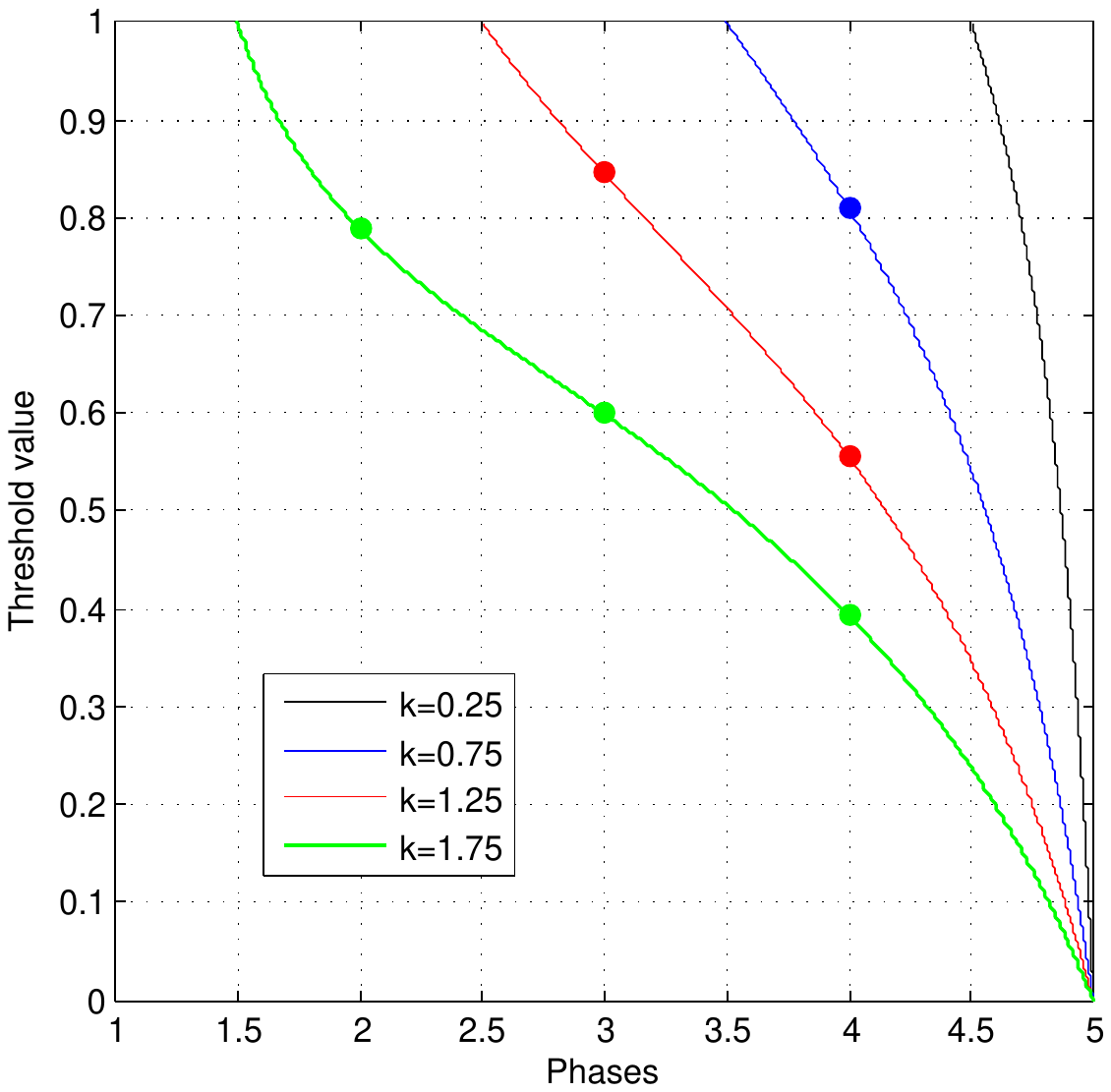}
\vspace{-6.5cm}
\end{center}
\caption{Threshold values \eqref{thre} (in big dots) for different values of $k$ for $N=5$ phases, $\alpha=\frac{1}{2},\beta=1$.}
\label{Compkph}
\end{figure}

Let us study closer this tendency behavior and the threshold values \eqref{thre} for an specific case. We fix $N=5$ phases and the parameters $\alpha=\frac{1}{2},\beta=1$. The interval $(0,2)$ in Figure \ref{Figkk} is split in four subintervals $(0,1/2]\cup(1/2,1]\cup(1,3/2]\cup(3/2,1)$. We study four values of $k$ taken in the middle of each subinterval. In Figure \ref{Compkph} we have plotted the threshold function \eqref{thre} (depending on the phase $i$) for these four values of $k$. We observe that the threshold value at the boundary phase 5 is 0, meaning that the next phase of the process is always 4, and the threshold value at the boundary phase 1 is $\infty$, meaning that the next phase of the process is always 2.

We observe, for $k\in(0,1/2)$ in Figure \ref{Compkph}, that there is no threshold between 0 and 1, i.e. always a forward tendency. In the following table we can confirm this behavior for a sample path of the process. Here in the first row we have written the sequence of transitions of phases and in the second row the position of $X_t$ at the end of that phase.

\medskip

\begin{tabular}{|c|c|c|c|c|c|c|c|c|c|c|c|c|c|c|c|c|c|c|c|}
  \hline
  Phases&1&2&3&2&3&4&5&4&5&4&5&4&5&4&5
 \\
  \hline
  Position &.68&.79&.82&.80&.69&.71&.76&.55&.86&.86&.85&.71&.59&.76&.69\\
  \hline
\end{tabular}

\medskip

For $k\in(1/2,1)$ in Figure \ref{Compkph}, only phase 4 has a threshold value (around 0.8). That means that for phases 2 and 3 there is always a forward tendency. At the moment of transition at the end of phase 4, if the position of $X_t$ is greater that the threshold value 0.8, then there is a backward tendency. Otherwise, there is a forward tendency. In the following table, as before, we have run a sample path for this situation:

\medskip

\begin{tabular}{|c|c|c|c|c|c|c|c|c|c|c|c|c|c|c|c|c|c|c|c|}
  \hline
Phases&  1&2&3&2&3&4&3&4&5&4&5&4&5&4&5
 \\
  \hline
 Position &.49&.53&.56&.51&.46&.87&.52&.45&.12&.13&.20&.04&.13&.12&.37\\
  \hline
\end{tabular}

\medskip

For $k\in(1,3/2)$ in Figure \ref{Compkph}, phases 3 and 4 have a threshold value (around 0.85 and 0.55 respectively). Again, phase 2 has a forward tendency, while phases 3 and 4 depend on the position of $X_t$ as before. In the following table we have run a sample path for this situation:
\medskip

\begin{tabular}{|c|c|c|c|c|c|c|c|c|c|c|c|c|c|c|c|c|c|c|c|}
  \hline
  Phases&1&2&3&4&5&4&3&2&3&2&3&4&5&4&3 \\
  \hline
  Position &.71&.83&.78&.43&.83&.81&.88&.78&.64&.86&.69&.72&.42&.58&.31\\
  \hline
\end{tabular}

\medskip

Finally, for $k\in(3/2,1)$ in Figure \ref{Compkph}, all middle phases have a threshold value (around 0.8, 0.6 and 0.4 respectively). This is the situation with the maximum backward tendency, but it depends strongly on the position of $X_t$. Again, in the following table we have run a sample path for this situation:

\medskip

\begin{tabular}{|c|c|c|c|c|c|c|c|c|c|c|c|c|c|c|c|c|c|c|c|}
  \hline
  Phases&1&2&1&2&3&2&3&4&3&2&1&2&3&4&3 \\
  \hline
  Position &.56&.72&.43&.47&.75&.66&.64&.87&.55&.62&.65&.17&.57&.62&.43 \\
  \hline
\end{tabular}

\bigskip

Finally, we will give the explicit expression of the vector-valued invariant distribution $\bm\psi(y)$ (see the end of Section \ref{NEW} for definitions). We can get an explicit formula in this case since $\mathcal{A}$, given by \eqref{DD}, is self-adjoint and we have an explicit expression of the orthogonality measure $\bm W$. Formula \eqref{IDW} gives
\begin{equation}\label{psijv}
    \bm\psi(y)=\left(\int_0^1\bm e^*_N\bm W(y)\bm e_Ndy\right)^{-1}\bm e^*_N\bm W(y).
\end{equation}
where $\bm W$ is given by \eqref{WW}. A closer look to the 0-norm and $\bm H$ shows, for $j=1,\ldots,N$,
\begin{equation}\label{psij}
\psi_j(y)=y^{\alpha+N-j}(1-y)^{\beta}\begin{pmatrix}
                                      N-1 \\
                                       j-1 \\
                                     \end{pmatrix}\begin{pmatrix}
                                       \alpha+\beta+N \\
                                       \alpha \\
                                     \end{pmatrix}\frac{(\beta+N)(k)_{N-j}(\beta-k+1)_{j-1}}{(\alpha+\beta-k+2)_{N-1}},
\end{equation}
where $(a)_n$ denotes the Pochhammer symbol defined by $(a)_n=a(a+1)\cdots(a+n-1)$ for $n>0$, $(a)_0=1$. Another possibility of obtaining \eqref{psijv} is by looking the spectral representation of $\bm P(t;x,y)$ \eqref{SRPDex1}. From this representation we observe that the limit as $t\rightarrow\infty$ only depends on the eigenvalues $\bm\Gamma_n$, given in \eqref{Lamb}. All diagonal entries of $\bm\Gamma_n$ are strictly negative for the range of all parameters $n, \alpha, \beta, k$ and $N$ involved, \emph{except} for $n=0$ and entry $(1,1)$ of $\bm\Gamma_0$, which is 0. Therefore taking limits in \eqref{SRPDex1} when $t\rightarrow\infty$, we have that the only relevant element is the first term of the sum, given by
$$
\lim_{t\rightarrow\infty}\bm P(t;x,y)=\bm\Phi_0(x)\bm E_{11}\bm\Phi_0^*(y)\bm W(y).
$$
An straightforward computation of the expression above, using the notation given in Section \ref{1step}, leads to a matrix where all rows are equal to \eqref{psijv}.

%We recall that $\bm E_{ij}$ denotes the matrix with 1 at entry $(i,j)$ and 0 otherwise. An straightforward computation of the expression above, using the notation given in Section \ref{1step}, leads to
%\begin{equation*}\label{InvDist}
%  \lim_{t\rightarrow\infty}\bm P(t;x,y)=\left(\|\bm R_0\|_{\widetilde{\bm W}}^{-2}\right)_{1,1}y^{\alpha}(1-y)^{\beta}\bm e_N\bm e_N^*\bm Hy^{\bm J}.
%\end{equation*}
%Therefore the invariant distribution is given by any row of the above expression, i.e.
%\begin{equation}\label{psijv}
%\bm\psi(y)=\left(\|\bm R_0\|_{\widetilde{\bm W}}^{-2}\right)_{1,1}y^{\alpha}(1-y)^{\beta}\bm e_N^*\bm Hy^{\bm J}.
%\end{equation}

This vector-valued invariant distribution \eqref{psijv} is valid only when the process is positive recurrent, i.e. $\alpha,\beta\geq0$. For $-1<\alpha,\beta<0$, \eqref{psijv} is also meaningful, but the boundary points of the process are now absorbing, meaning that the correct vector-valued invariant distribution of such cases involves mass jumps at the boundaries 0 and 1 plus a density portion of the form \eqref{psijv}.

\begin{figure}[h]
\begin{center}
\vspace{-5.0cm}
\includegraphics[height=20cm]{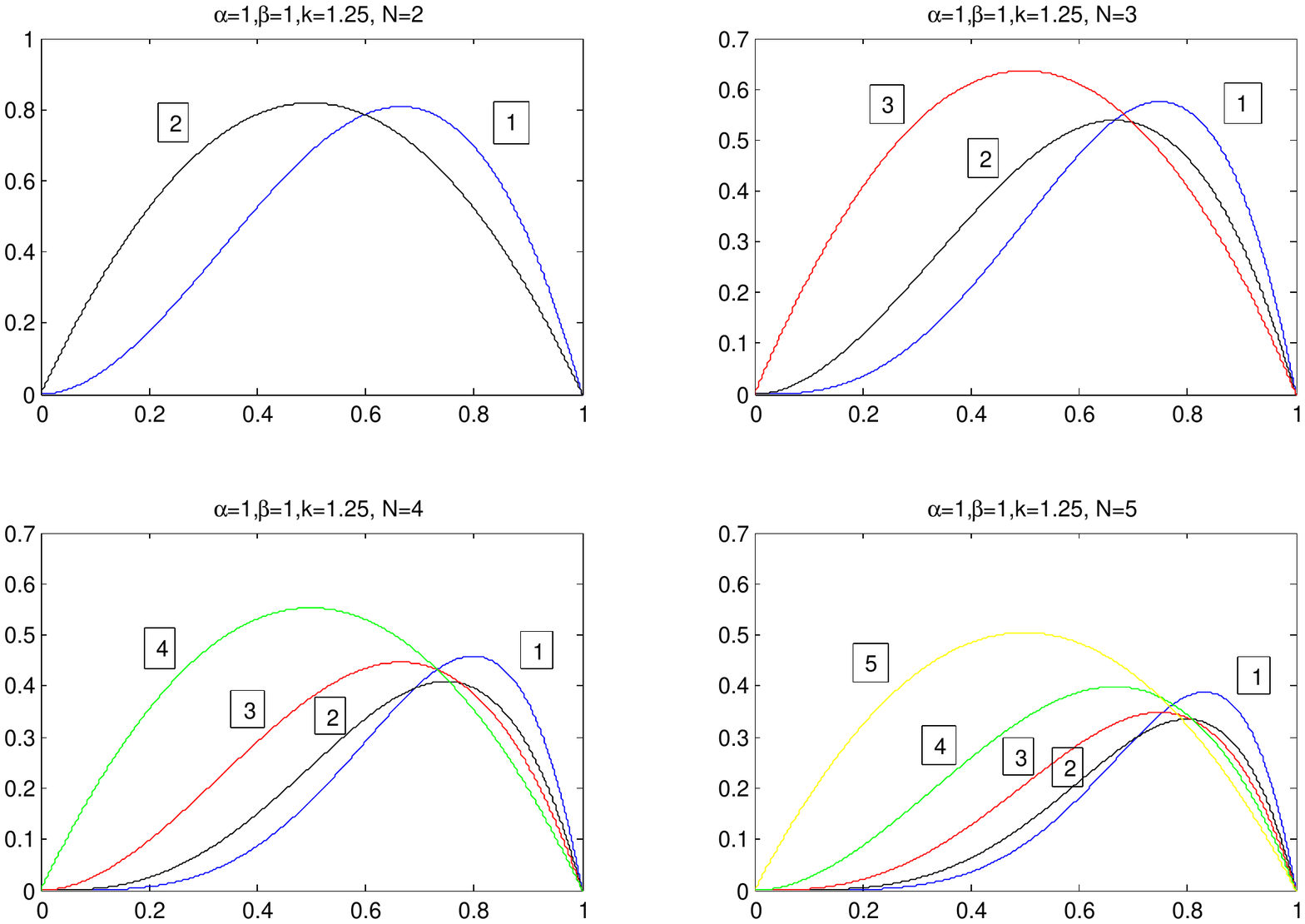}
\vspace{-6.0cm}
\end{center}
\caption{The components of the vector-valued invariant distribution $\psi_j(y), j=1,\ldots,N$ (see \eqref{psij}), for $N=2,3,4,5,$ and $\alpha, \beta=1, k=5/4.$}
\label{IDs}
\end{figure}

As an illustration, in Figure \ref{IDs}, we have plotted $\psi_j(y)$, $j=1,\ldots N$, in \eqref{psij} for $N=2,3,4,5$. We have fixed the parameters to $\alpha, \beta=1, k=5/4$. We remark that $\int_0^1\psi_j(y)dy<1$ for each component, but $\sum_{j=1}^N\int_0^1\psi_j(y)dy=1$.

From these plots we immediately see that, for a large time, if the process is at one state in $(0,1/2)$, then there is a very high probability that the process is located at one of the last phases. Nevertheless if the process is at one state in $(1/2,1)$, there is more or less the same probability that the process is located at any phase (higher probability for the first phases).

This interpretation may change depending on the value of $k$. If $k$ is close to $\beta+1$ then we have the maximum backward tendency and the process is more likely to be at one of the initial phases. Nevertheless, if $k$ is close to 0, then we have maximum forward tendency, and the process is more likely to be located at one of the final phases.

\subsection{Final remarks}\label{FINREM}

From the considerations in Section \ref{PROB} we are ready to give an interpretation of this bivariate Markov process in relation with the Wright-Fisher diffusion model involving only mutation effects. Consider a big population of constant size composed of two types A and B, where the intensity of mutation from A to B is given in terms of the parameter $\beta$ and the intensity of mutation from B to A is given in terms of the parameter $\alpha$. The value that the process $X_t$ takes in the interval $(0,1)$ is the rate $\#$ A$/(\#$ A+$\#$ B), i.e. if $X_t<1/2$ then $\#$ B$>\#$ A and if $X_t>1/2$ then $\#$ A$>\#$ B. We assume that the intensity of mutations $\alpha$ and $\beta,$ are more or less the same. If $\alpha>\beta$ then the population of A's tends to survive against the population of B's (the opposite for $\alpha<\beta$).

Now the process can move through $N$ different phases. At phase $N$ the process behaves as a regular Wright-Fisher model (in terms of the Jacobi diffusion). At phase $N-1$ the intensity $B\rightarrow A$ increases in one unity, at phase $N-2$ the intensity $B\rightarrow A$ increases in two unities, and so on, until the first phase, where the intensity $B\rightarrow A$ increases in $N-1$ unities. Therefore moving through the first phases means that the intensity of mutation $B\rightarrow A$ is much higher if we compare with the intensity of mutation $A\rightarrow B$, i.e. it helps to increase the population of A's against B's.

The meaning of the indicator parameter $k$ (depending on $\beta$) is the following. If $k$ is close to 0, then we have a forward tendency, i.e. the process tends to spend more time at the last phases, i.e. the populations of A's and B's ``fight'' more or less in the same conditions (remember that we are assuming that $\alpha$ is more or less equal to $\beta$). On the contrary, if this indicator $k$ is close to $\beta+1$, then we have a backward tendency, i.e. the process tends to spend more time at the initial phases, i.e. it helps the population of A's to survive against the population of B's. For the middle values of $k$ there can be forward/backward tendency, and it depends on the number of A's and B's, i.e. the position of $X_t$.

\section{Conclusions}

In this paper we have studied bivariate Markov processes with diffusion and discrete components from its matrix-valued infinitesimal operator properties. This infinitesimal operator is a second-order differential operator with matrix-valued coefficients with certain properties. We have given some theoretical results related with these processes, like the backward/forward equations, the spectral representation in terms of the eigenfunctions of the infinitesimal operator, recurrence considerations and the invariant distribution. All these results have been applied to an example coming from group representation theory and related with matrix-valued spherical functions, which can be viewed as a variant of the Wright-Fisher diffusion model involving only mutation effects. Since we have an explicit set of eigenfunctions with the corresponding orthogonality measure we can approximate very accurately the matrix-valued probability density $\bm P(t;x,y)$ and get an explicit expression of the vector-valued invariant distribution, among other properties.

Certainly, new examples of bivariate Markov processes may be possible, specially when the diffusion state space is other than a bounded interval. In \cite{GPT1} one can find other types of examples. The example studied in this paper, also known as the \emph{one-step example}, is just one class of a very large class of matrix-valued spherical functions associated with the complex projective space. The second-order differential operators associated with these matrix-valued spherical functions seem to fit very well in the theory of bivariate Markov processes. For instance, in \cite{GPT6}, one can find a case of the \emph{two-steps example} ($4\times4$). The second-order differential operator is also the infinitesimal operator of a bivariate Markov process. In this case, apart from the parameters $\alpha, \beta>-1$, there are two new extra parameters, $k_1$ and $k_2$ with the condition that $0<k_1<k_2<\beta+1$. Therefore this example can be viewed as another variant of the Wright-Fisher model involving only mutation effects, but with two parameters to study, which decide how the process move through all the phases.

Eigenfunctions of second-order differential operators with matrix-valued coefficients have been studied in the context of matrix-valued orthogonal polynomials, see \cite{DG1, GPT5, GPT6}. These examples are potential candidates of bivariate Markov processes. The difference is that these differential operators do not exactly fit into the framework of bivariate Markov processes. An inverse transformation of the type \eqref{Relat} is going to be needed in order to get an infinitesimal operator of a bivariate Markov process, which in principle is not an easy task. This consideration, however, goes beyond the scope of this paper, and will be pursued elsewhere.


\begin{thebibliography}{sssssss}

\bibitem{AAR}\textrm{Andrews, G. E., Askey, R. and Roy, R.},
\textit{Special functions}, Encyclopedia of Mathematics, No. 71.
Cambridge, (1999).

\bibitem{B} Berman, S. M., {\em
A bivariate Markov process with diffusion and discrete components}, Comm. Statist. Stochastic Models \textbf{10} (1994), No. 2, 271--308.

\bibitem{BW} Bhattacharya, R. N. and Waymire, E. C., {\em
Stochastic processes with applications}, Wiley Series in Probability and Mathematical Statistics: Applied Probability and Statistics. John Wiley \& Sons, Inc., New York, 1990.

%\bibitem{Bl} Blackwell, P. G., {\em
%Bayesian inference for Markov processes with diffusion and discrete components}, Biometrika \textbf{90} (2003), No. 3, 613--627.

\bibitem{DPS} Damanik, D., Pushnitski, A. and Simon, B., \emph{The analytic theory of matrix
orthogonal polynomials}, Surv. Approx. Theory $\mathbf{4}$ (2008), 1--85.

%\bibitem{DRSZ} Dette, H., Reuther, B., Studden, W. and Zygmunt, M., {\em
%Matrix measures and random walks with a block tridiagonal transition
%matrix}, SIAM J. Matrix Anal. Applic. \textbf{29}, No. 1 (2006) pp.
%117--142.

%\bibitem{DR} Dette, H. and Reuther, B., {\em
%Some comments on quasi-birth-and-death processes and matrix measures}, J. Probability and Statistics
%Volume 2010 (2010), Article ID 730543, 23 pages.

\bibitem{DG1}\textrm{Dur\'{a}n, A. J. and Gr\"{u}nbaum, F. A.},
\textit{Orthogonal matrix polynomials satisfying second order
differential equations}, Internat. Math. Research Notices, 2004:
$\mathbf{10}$ (2004), 461--484.

\bibitem{GHO}\textrm{Ghosh, M. K., Arapostathis, A. and Marcus, S. I.},
\textit{Ergodic control of switching diffusions}, SIAM J. Control Optim.,
$\mathbf{35}$ (1997), No. 6, 1952--1988.

\bibitem{GH1}\textrm{Griego, R. J. and Hersh, R.},
\textit{Random evolutions, Markov chains, and systems of
partial differential equations}, Proc. Nat. Acad. Sci. U.S.A.,
$\mathbf{62}$ (1969), 305--308.


%\bibitem{G2} Gr\"unbaum, F. A., {\em
%Random walks and orthogonal polynomials: some challenges},
%in Probability, Geometry and Integrable Systems, Mark Pinsky and Bjorn Birnir editors, MSRI Publication,
%vol. \textbf{55}, 2007, 241--260. See also arXiv:math/0703375v1.

\bibitem{GdI2} Gr\"unbaum, F. A. and de la Iglesia, M. D.,
\textit{Matrix valued orthogonal polynomials arising from group
representation theory and a family of quasi-birth-and-death
processes},  SIAM J. Matrix Anal. Appl.  \textbf{30}  (2008),  No. 2, 741--761.

\bibitem{GPT1} Gr\"unbaum, F. A., Pacharoni, I., and Tirao, J. A., {\em
Matrix valued spherical functions associated to the complex
projective plane}, J. Functional Analysis {\bf 188} (2002), 350--441.

\bibitem{GPT5} Gr\"unbaum, F. A., Pacharoni, I. and Tirao, J. A.,
 {\em Matrix valued orthogonal polynomials of the Jacobi
type}, Indag. Mathem. {\bf 14} (2003), Nos. 3,4, 353--366.

\bibitem{GPT6} Gr\"unbaum, F. A., Pacharoni, I. and Tirao, J.
A., {\em Matrix valued orthogonal polynomials of the Jacobi type:
the role of group representation theory}, Ann. Inst. Fourier, {\bf
55} (2005), No. 5, 1--18.

\bibitem{GPT7} Gr\"unbaum, F. A., Pacharoni, I. and Tirao, J.
A., {\em Two stochastic models of a random walk in the $U(n)$-spherical duals of $U(n + 1)$}, see arXiv:1010.0720 (2010).

\bibitem{He1} Hersh, R., {\em Random evolutions: a survey of results and problems},  Rocky
Mountain J. Math., {\bf 4} (1974), 443--477.

\bibitem{He2} Hersh, R., {\em The birth of random evolutions},  Mathematical Intelligencer, {\bf 25} (2003), 53--60.

\bibitem{KT} Karlin, S. and Taylor, H. M., {\em
A second course in stochastic processes}, Academic Press, Inc. Harcourt Brace Jovanovich, Publishers, New York-London, 1981.

%\bibitem{LR1} Latouche, G. and Ramaswami, V., {\em Introduction to matrix
%analytic methods in stochastic modeling}, ASA-SIAM Series on
%Statistics and Applied Probability, 1999.

%\bibitem{N} Neuts, M. F., {\em Structured stochastic matrices of $M/G/1$ type and
%their applications}, Marcel Dekker, NY, 1989.

\bibitem{PT1}  Pacharoni, I. and  Tirao, J. A.,
{\em Matrix valued orthogonal polynomials arising from the complex
projective space}, Constr. Approx.  \textbf{25} (2007), 177--192.

\bibitem{Pap} Papanicolaou, G. C., {\em Random media}, Springer-Verlag, Berlin, 1987.

\bibitem{Pins} Pinsky, M. A., {\em Lectures on random evolutions}, World Scientific, Singapore, 1991.

\bibitem{Sch} Schoutens, W., {\em Stochastic processes and orthogonal polynomials}, Lecture Notes in Statistics, 146.
Springer-Verlag, New York, 2000.

\bibitem{Swi} Swishchuk, A., {\em Random evolutions and their applications. New Trends}, Kluwer AP, Dordrecht, 2000.


%\bibitem{VW}  Varadarajan, V. S. and  Weisbart, D.,
%{\em Convergence of quantum systems on grids}, J. Math. Anal. Appl. \textbf{336} (2007), pp.
%608--624.


\end{thebibliography}
\end{document}